%
% The projective geometry of the  Gale Transform
%
% David Eisenbud, Sorin Popescu
%
% (version of 07/22/98)
%
% Plain TeX with  diagrams.tex

%-----------------------------------------------
% begin.tex
%-----------------------------------------------
%=========================================================================%
% load begin.tex only once, but keep count to match \bye commands
%=========================================================================%

\ifx\begin\undefined\else\global\advance\srcdepth by
1\expandafter\endinput\fi

\def\begin{}
\newcount\srcdepth
\srcdepth=1

\outer\def\bye{\global\advance\srcdepth by -1
  \ifnum\srcdepth=0
    \def\endcmd{\vfill\eject\nopagenumbers\par\vfill\supereject\end}
  \else\def\endcmd{}\fi
  \endcmd
}

%=========================================================================%
% initialize TeX
%=========================================================================%

\magnification=\magstephalf
\baselineskip=13pt
\hsize = 5.5truein
\hoffset = 0.5truein
\vsize = 8.5truein
\voffset = 0.2truein
\emergencystretch = 0.05\hsize

\newif\ifblackboardbold

% comment out the following line if AMS msbm fonts aren't available
\blackboardboldtrue

%=========================================================================%
% select fonts
%=========================================================================%

\font\sectionfont=cmbx12

% Establish AMS blackboard bold fonts without using amssym.def, amssym.tex

\newfam\bboldfam
\ifblackboardbold
\font\tenbbold=msbm10
\font\sevenbbold=msbm7
\font\fivebbold=msbm5
\textfont\bboldfam=\tenbbold
\scriptfont\bboldfam=\sevenbbold
\scriptscriptfont\bboldfam=\fivebbold
\def\bbold{\fam\bboldfam\tenbbold}
\else
\def\bbold{\bf}
\fi

%=========================================================================%
% font size-changing command ("A Beginner's Book of TeX" p35, p275)
%=========================================================================%

\font\Arm=cmr8
\font\Ai=cmmi8
\font\Asy=cmsy8
\font\Abf=cmbx8
\font\Brm=cmr6
\font\Bi=cmmi6
\font\Bsy=cmsy6
\font\Bbf=cmbx6
\font\Crm=cmr5
\font\Ci=cmmi5
\font\Csy=cmsy5
\font\Cbf=cmbx5

\ifblackboardbold
\font\Abbold=msbm10 at 8pt
\font\Bbbold=msbm7 at 6pt
\font\Cbbold=msbm5
\fi

\def\smallmath{%
\textfont0=\Arm \scriptfont0=\Brm \scriptscriptfont0=\Crm
\textfont1=\Ai \scriptfont1=\Bi \scriptscriptfont1=\Ci
\textfont2=\Asy \scriptfont2=\Bsy \scriptscriptfont2=\Csy
\textfont\bffam=\Abf \scriptfont\bffam=\Bbf \scriptscriptfont\bffam=\Cbf
\def\rm{\fam0\Arm}\def\mit{\fam1}\def\oldstyle{\fam1\Ai}%
\def\bf{\fam\bffam\Abf}%
\ifblackboardbold
\textfont\bboldfam=\Abbold
\scriptfont\bboldfam=\Bbbold
\scriptscriptfont\bboldfam=\Cbbold
\def\bbold{\fam\bboldfam\Abbold}%
\fi
}

%=========================================================================%
% single-pass symbolic theorem labeling
%=========================================================================%

% Because this is a single-pass mechanism with no .aux file, forward
% references need to be declared in advance:

%   \forward{thm:main}{Theorem}{1.1}

% This is also the mechanism for "timely" declaration of labels, which
% will usually be buried within the corresponding theorem macros.
% A warning is issued if a label redeclaration is inconsistent, allowing
% forward references to be manually fixed.

%   \ref{thm:main} produces "Theorem~1.1"
%   \refs{thm:main} produces "Theorems~1.1"
%   \refn{thm:main} produces "1.1"

% Some TeX adapted from "The Advanced TeXbook" by David Salomon, chapter 9.

% Implementers: The code for \forward is subtle. Its second argument must
% be provided literally, e.g. "Theorem" rather that "\capitalize{theorem}".
% Its third argument must either be literal or a macro that expands
% directly to a literal, e.g. "\edef\numtoks{\number\proccount}".
% This use of \edef cannot be replaced by \def, which defers expansion.
% Failure to follow these rules will cause spurious warnings that forward
% references are inconsistent, when they are in fact consistent after
% expansion. Note the "Towers of Palo Alto" recreational math problem
% involving the iterated use of \expandafter to expand the first argument
% to \forwardsub before calling it.

\newlinechar=`@
\def\forwardmsg#1#2#3{\immediate\write16{@*!*!*!* forward reference should
be: @\noexpand\forward{#1}{#2}{#3}@}}
\def\nodefmsg#1{\immediate\write16{@*!*!*!* #1 is an undefined reference@}}

\def\forwardsub#1#2{\def\newref{{#2}{#1}}}

\def\forward#1#2#3{%
\expandafter\expandafter\expandafter\forwardsub\expandafter{#3}{#2}
\expandafter\ifx\csname#1\endcsname\relax\else%
\expandafter\ifx\csname#1\endcsname\newref\else%
\forwardmsg{#1}{#2}{#3}\fi\fi%
\expandafter\let\csname#1\endcsname\newref}

\def\firstarg#1{\expandafter\argone #1}\def\argone#1#2{#1}
\def\secondarg#1{\expandafter\argtwo #1}\def\argtwo#1#2{#2}

\def\ref#1{\expandafter\ifx\csname#1\endcsname\relax
  {\nodefmsg{#1}\bf`#1'}\else
  \expandafter\firstarg\csname#1\endcsname
  ~\expandafter\secondarg\csname#1\endcsname\fi}

\def\refs#1{\expandafter\ifx\csname#1\endcsname\relax
  {\nodefmsg{#1}\bf`#1'}\else
  \expandafter\firstarg\csname #1\endcsname
  s~\expandafter\secondarg\csname#1\endcsname\fi}

\def\refn#1{\expandafter\ifx\csname#1\endcsname\relax
  {\nodefmsg{#1}\bf`#1'}\else
  \expandafter\secondarg\csname #1\endcsname\fi}

%=========================================================================%
% widow control
%=========================================================================%

% usage:
% \widow{.2} % start new page if <.2 page left

\def\widow#1{\vskip 0pt plus#1\vsize\goodbreak\vskip 0pt plus-#1\vsize}

%=========================================================================%
% sections and theorems
%=========================================================================%

% use \showlabels or \showlabelsabove to display section and theorem labels

\def\marginlabel#1{}

\def\showlabelsabove{
\font\labelfont=cmss10 at 6pt
\def\marginlabel##1{\rlap{\smash{\raise 10pt\hbox{\labelfont##1}}}}
}

\newcount\seccount
\newcount\proccount
\seccount=0
\proccount=0

\def\stdskip{\vskip 9pt plus3pt minus 3pt}
\def\stdbreak{\par\removelastskip\penalty-100\stdskip}

\def\proof{\stdbreak\noindent{\sl Proof. }}

\def\qed{\vrule height 1.2ex width .9ex depth .1ex}

\def\Box{
  \ifmmode\eqno\qed
  \else\ifvmode\removelastskip\line{\hfil\qed}
  \else\unskip\quad\hskip-\hsize
    \hbox{}\hskip\hsize minus 1em\qed\par
  \fi\stdbreak\fi}

\def\ifempty#1#2\endB{\ifx#1\endA}
\def\makeref#1#2#3{\ifempty#1\endA\endB\else\forward{#1}{#2}{#3}\fi}

\outer\def\section#1 #2\par{
  \removelastskip
  \global\advance\seccount by 1
  \global\proccount=0\relax
                \edef\numtoks{\number\seccount}
  \makeref{#1}{Section}{\numtoks}
  \widow{.05}
  \vskip 24pt plus 6pt minus 6 pt
  \message{#2}
  \leftline{\marginlabel{#1}\sectionfont\numtoks\quad #2}
  \nobreak\stdskip}

\def\proclamation#1#2{
  \outer\expandafter\def\csname#1\endcsname##1 ##2\par{
  \stdbreak
  \advance\proccount by 1
  \edef\numtoks{\number\seccount.\number\proccount}
  \makeref{##1}{#2}{\numtoks}
  \noindent{\marginlabel{##1}\bf #2 \numtoks\enspace}
  {\sl##2\par}
  \stdbreak}}

\def\othernumbered#1#2{
  \outer\expandafter\def\csname#1\endcsname##1{
  \stdbreak
  \advance\proccount by 1
  \edef\numtoks{\number\seccount.\number\proccount}
  \makeref{##1}{#2}{\numtoks}
  \noindent{\marginlabel{##1}\bf #2 \numtoks\enspace}}}

\proclamation{definition}{Definition}
\proclamation{lemma}{Lemma}
\proclamation{proposition}{Proposition}
\proclamation{theorem}{Theorem}
\proclamation{corollary}{Corollary}
\proclamation{conjecture}{Conjecture}

\othernumbered{example}{Example}
\othernumbered{remark}{Remark}
\othernumbered{construction}{Construction}
\othernumbered{problem}{Problem}
%=========================================================================%
% enable postscript illustrations using epsf.tex
%=========================================================================%

% Usage:
% \draw{70}{fig}{} % draw fig.eps at 70% scale
% \draw{999}{fig}{} % draw fig.eps scaled to width of page

% Optional third argument can be multiple calls to \figtext; see below.
% More generally, the third argument is read in vertical mode, with the
% reference point at the lower left corner of the eps picture, whose
% dimensions are contained in the dimen registers \drawx and \drawy.
% This enables using TeX to generate the text that goes with the picture.
% To request that the picture be widened to respect the added text,
% examine and modify the dimen registers \ngap, \egap, \sgap, \wgap.
% This is done automatically by the \figtext macro.

% These macros rely on "epsf.tex" which is the lowest level interface
% available for including encapsulated Postscript files in TeX documents.
% Rather that manually reading the .eps file to compute the nominal size,
% the \epsfbox macro is called twice, and two of its internal registers
% are examined after the first call. A major change to epsf.tex (unlikely)
% will require changes here.

\input epsf

\newcount\figcount
\figcount=0
\newbox\drawing
\newcount\drawbp
\newdimen\drawx
\newdimen\drawy
\newdimen\ngap
\newdimen\sgap
\newdimen\wgap
\newdimen\egap

\def\drawbox#1#2#3{\vbox{
  \setbox\drawing=\vbox{\offinterlineskip\epsfbox{#2.eps}\kern 0pt}
  \drawbp=\epsfurx
  \advance\drawbp by-\epsfllx\relax
  \multiply\drawbp by #1
  \divide\drawbp by 100
  \drawx=\drawbp truebp
  \ifdim\drawx>\hsize\drawx=\hsize\fi
  \epsfxsize=\drawx
  \setbox\drawing=\vbox{\offinterlineskip\epsfbox{#2.eps}\kern 0pt}
  \drawx=\wd\drawing
  \drawy=\ht\drawing
  \ngap=0pt \sgap=0pt \wgap=0pt \egap=0pt
  \setbox0=\vbox{\offinterlineskip
    \box\drawing \ifgridlines\drawgrid\drawx\drawy\fi #3}
  \kern\ngap\hbox{\kern\wgap\box0\kern\egap}\kern\sgap}}

\def\draw#1#2#3{
  \setbox\drawing=\drawbox{#1}{#2}{#3}
  \advance\figcount by 1
  \goodbreak
  \midinsert
  \centerline{\ifgridlines\boxgrid\drawing\fi\box\drawing}
  \smallskip
  \vbox{\offinterlineskip
    \centerline{Figure~\number\figcount}
    \smash{\marginlabel{#2}}}
  \endinsert}

\def\nextfigtoks{%
  \advance\figcount by 1%
  \edef\numtoks{\number\figcount}%
  \advance\figcount by -1}

\newif\ifgridlines
\newbox\figtbox
\newbox\figgbox
\newdimen\figtx
\newdimen\figty

\newdimen\bwd
\bwd=2sp % 2sp (1/32768") is smallest visible width for Textures

\def\hline#1{\vbox{\smash{\hbox to #1{\leaders\hrule height \bwd\hfil}}}}

\def\vline#1{\hbox to 0pt{%
  \hss\vbox to #1{\leaders\vrule width \bwd\vfil}\hss}}

\def\clap#1{\hbox to 0pt{\hss#1\hss}}
\def\vclap#1{\vbox to 0pt{\offinterlineskip\vss#1\vss}}

\def\hstutter#1#2{\hbox{%
  \setbox0=\hbox{#1}%
  \hbox to #2\wd0{\leaders\box0\hfil}}}

\def\vstutter#1#2{\vbox{
  \setbox0=\vbox{\offinterlineskip #1}
  \dp0=0pt
  \vbox to #2\ht0{\leaders\box0\vfil}}}

\def\crosshairs#1#2{
  \dimen1=.002\drawx
  \dimen2=.002\drawy
  \ifdim\dimen1<\dimen2\dimen3\dimen1\else\dimen3\dimen2\fi
  \setbox1=\vclap{\vline{2\dimen3}}
  \setbox2=\clap{\hline{2\dimen3}}
  \setbox3=\hstutter{\kern\dimen1\box1}{4}
  \setbox4=\vstutter{\kern\dimen2\box2}{4}
  \setbox1=\vclap{\vline{4\dimen3}}
  \setbox2=\clap{\hline{4\dimen3}}
  \setbox5=\clap{\copy1\hstutter{\box3\kern\dimen1\box1}{6}}
  \setbox6=\vclap{\copy2\vstutter{\box4\kern\dimen2\box2}{6}}
  \setbox1=\vbox{\offinterlineskip\box5\box6}
  \smash{\vbox to #2{\hbox to #1{\hss\box1}\vss}}}

\def\boxgrid#1{\rlap{\vbox{\offinterlineskip
  \setbox0=\hline{\wd#1}
  \setbox1=\vline{\ht#1}
  \smash{\vbox to \ht#1{\offinterlineskip\copy0\vfil\box0}}
  \smash{\vbox{\hbox to \wd#1{\copy1\hfil\box1}}}}}}

\def\drawgrid#1#2{\vbox{\offinterlineskip
  \dimen0=\drawx
  \dimen1=\drawy
  \divide\dimen0 by 10
  \divide\dimen1 by 10
  \setbox0=\hline\drawx
  \setbox1=\vline\drawy
  \smash{\vbox{\offinterlineskip
    \copy0\vstutter{\kern\dimen1\box0}{10}}}
  \smash{\hbox{\copy1\hstutter{\kern\dimen0\box1}{10}}}}}

\def\figtext#1#2#3#4#5{
  \setbox\figtbox=\hbox{#5}
  \dp\figtbox=0pt
  \figtx=-#3\wd\figtbox \figty=-#4\ht\figtbox
  \advance\figtx by #1\drawx \advance\figty by #2\drawy
  \dimen0=\figtx \advance\dimen0 by\wd\figtbox \advance\dimen0 by-\drawx
  \ifdim\dimen0>\egap\global\egap=\dimen0\fi
  \dimen0=\figty \advance\dimen0 by\ht\figtbox \advance\dimen0 by-\drawy
  \ifdim\dimen0>\ngap\global\ngap=\dimen0\fi
  \dimen0=-\figtx
  \ifdim\dimen0>\wgap\global\wgap=\dimen0\fi
  \dimen0=-\figty
  \ifdim\dimen0>\sgap\global\sgap=\dimen0\fi
  \smash{\rlap{\vbox{\offinterlineskip
    \hbox{\hbox to \figtx{}\ifgridlines\boxgrid\figtbox\fi\box\figtbox}
    \vbox to \figty{}
    \ifgridlines\crosshairs{#1\drawx}{#2\drawy}\fi
    \kern 0pt}}}}

% macros to add space to text on specified sides

\def\hpad#1#2#3{\hbox{\kern #1\hbox{#3}\kern #2}}
\def\vpad#1#2#3{\setbox0=\hbox{#3}\dp0=0pt\vbox{\kern #1\box0\kern #2}}

% macro to give one text string the apparent height of another

% macro to center one text string over another

\def\stack#1#2#3{\vbox{\offinterlineskip
  \setbox2=\hbox{#2}
  \setbox3=\hbox{#3}
  \dimen0=\ifdim\wd2>\wd3\wd2\else\wd3\fi
  \hbox to \dimen0{\hss\box2\hss}
  \kern #1
  \hbox to \dimen0{\hss\box3\hss}}}

% macros to hide size of trailing exponents

\def\hexp#1{%
  \setbox0=\hbox{${}^{#1}$}%
  \hbox to .5\wd0{\box0\hss}}

%=========================================================================%
% macros for matrices and arrows
%=========================================================================%

% typical usage:
%   \rightarrowmat{2pt}{4pt}{d & bd \cr \!-c & 0 \cr 0 & -ac \cr}

\def\bmatrix#1#2{{\smallmath\left[\vcenter{\halign
  {&\kern#1\hfil$##\mathstrut$\kern#1\cr#2}}\right]}}

\def\rightarrowmat#1#2#3{
  \setbox1=\hbox{\kern#2$\bmatrix{#1}{#3}$\kern#2}
  \,\vbox{\offinterlineskip\hbox to\wd1{\hfil\copy1\hfil}
    \kern 3pt\hbox to\wd1{\rightarrowfill}}\,}

\def\leftarrowmat#1#2#3{
  \setbox1=\hbox{\kern#2$\bmatrix{#1}{#3}$\kern#2}
  \,\vbox{\offinterlineskip\hbox to\wd1{\hfil\copy1\hfil}
    \kern 3pt\hbox to\wd1{\leftarrowfill}}\,}

\def\rightarrowbox#1#2{
  \setbox1=\hbox{\kern#1\hbox{\smallmath #2}\kern#1}
  \,\vbox{\offinterlineskip\hbox to\wd1{\hfil\copy1\hfil}
    \kern 3pt\hbox to\wd1{\rightarrowfill}}\,}

\def\leftarrowbox#1#2{
  \setbox1=\hbox{\kern#1\hbox{\smallmath #2}\kern#1}
  \,\vbox{\offinterlineskip\hbox to\wd1{\hfil\copy1\hfil}
    \kern 3pt\hbox to\wd1{\leftarrowfill}}\,}

%=========================================================================%
% quire macros for preview mode and making booklets
%=========================================================================%

% \legalbooklet{20} makes a booklet from legal paper in landscape
% orientation, where "20" is the page count. To preview, give a negative
% pagecount. Either print using the legal duplex option on a modern laser
% printer, or struggle to simulate this effect manually. Bind using a long
% reach stapler.

% \preview squeezes two pages side by side in landscape orientation. It
% is not suitable for printing, but ideal for previewing on a two page
% monitor.

% \twoup squeezes two pages onto letter paper in landscape mode,
% suitable for printing.

% Each of these macros calls the file "quire.tex"

\def\bookletdims{
  \hsize=5.25truein
  \vsize=7truein
}

\def\legalbooklet#1{
  \input quire
  \bookletdims
  \htotal=7.0truein
  \vtotal=8.5truein
  % below computed from above
  \hoffset=\htotal
  \advance\hoffset by -\hsize
  \divide\hoffset by 2
  \voffset=\vtotal
  \advance\voffset by -\vsize
  \divide\voffset by 2
  \advance\voffset by -.0625truein
  \shhtotal=2\htotal
  % below doesn't need to change
  \horigin=0.0truein
  \vorigin=0.0truein
  \shstaplewidth=0.01pt
  \shstaplelength=0.66truein
  \shthickness=0pt
  \shoutline=0pt
  \shcrop=0pt
  \shvoffset=-1.0truein
  \ifnum#1>0\quire{#1}\else\qtwopages\fi
}

\def\preview{
  \input quire
  \bookletdims
  \hoffset=0.1truein
  \vtotal=8.5truein
  \shhtotal=14truein
  % below computed from above
  \voffset=\vtotal
  \advance\voffset by -\vsize
  \divide\voffset by 2
  \advance\voffset by -.0625truein
  \htotal=2\hoffset
  \advance\htotal by \hsize
  % below doesn't need to change
  \horigin=0.0truein
  \vorigin=0.0truein
  \shstaplewidth=0.5pt
  \shstaplelength=0.5\vtotal
  \shthickness=0pt
  \shoutline=0pt
  \shcrop=0pt
  \shvoffset=-1.0truein
  \qtwopages
}

\def\twoup{
  \input quire
  \hsize=4.79452truein % 5.25/1.095
  \vsize=7truein
  \vtotal=8.5truein
  \shhtotal=11truein
  % below computed from above
  \hoffset=-2\hsize
  \advance\hoffset by \shhtotal
  \divide\hoffset by 6
  \voffset=\vtotal
  \advance\voffset by -\vsize
  \divide\voffset by 2
  \advance\voffset by -12truept
  \htotal=2\hoffset
  \advance\htotal by \hsize
  % below doesn't need to change
  \horigin=0.0truein
  \vorigin=0.0truein
  \shstaplewidth=0.01pt
  \shstaplelength=0pt
  \shthickness=0pt
  \shoutline=0pt
  \shcrop=0pt
  \shvoffset=-1.0truein
  \qtwopages
}

%=========================================================================%
% timestamp (adapted from eplain.tex)
%=========================================================================%

\newcount\countA
\newcount\countB
\newcount\countC

\def\monthname{\begingroup
  \ifcase\number\month
    \or January\or February\or March\or April\or May\or June\or
    July\or August\or September\or October\or November\or December\fi
\endgroup}

\def\dayname{\begingroup
  \countA=\number\day
  \countB=\number\year
  \advance\countA by 0 % adjust after each leap day
  \advance\countA by \ifcase\month\or
    0\or 31\or 59\or 90\or 120\or 151\or
    181\or 212\or 243\or 273\or 304\or 334\fi
  \advance\countB by -1995
  \multiply\countB by 365
  \advance\countA by \countB
  \countB=\countA
  \divide\countB by 7
  \multiply\countB by 7
  \advance\countA by -\countB
  \advance\countA by 1
  \ifcase\countA\or Sunday\or Monday\or Tuesday\or Wednesday\or
    Thursday\or Friday\or Saturday\fi
\endgroup}

\def\timename{\begingroup
   \countA = \time
   \divide\countA by 60
   \countB = \countA
   \countC = \time
   \multiply\countA by 60
   \advance\countC by -\countA
   \ifnum\countC<10\toks1={0}\else\toks1={}\fi
   \ifnum\countB<12 \toks0={\sevenrm AM}
     \else\toks0={\sevenrm PM}\advance\countB by -12\fi
   \relax\ifnum\countB=0\countB=12\fi
   \hbox{\the\countB:\the\toks1 \the\countC \thinspace \the\toks0}
\endgroup}

\def\timestamp{\dayname, \the\day\ \monthname\ \the\year, \timename}

%==========================================================================
% macros (specific to this paper)
%==========================================================================

% surround with $ $ if not already in math mode
\def\enma#1{{\ifmmode#1\else$#1$\fi}}

% Gothic fonts from AMSTeX
\font\tengoth=eufm10  \font\fivegoth=eufm5
\font\sevengoth=eufm7
\newfam\gothfam  \scriptscriptfont\gothfam=\fivegoth
\textfont\gothfam=\tengoth \scriptfont\gothfam=\sevengoth
\def\goth{\fam\gothfam\tengoth}
%
% Bold italic fonts
\font\tenbi=cmmib10  \font\fivebi=cmmib5
\font\sevenbi=cmmib7
\newfam\bifam  \scriptscriptfont\bifam=\fivebi
\textfont\bifam=\tenbi \scriptfont\bifam=\sevenbi

\font\hd=cmbx10 scaled\magstep1
\font\abst=cmr9
\def \Box {\hfill\hbox{}\nobreak \vrule width 1.6mm height 1.6mm
depth 0mm  \par \goodbreak \smallskip}
\def \coker {\mathop{\rm coker}\nolimits}
\def \ker {\mathop{\rm ker}\nolimits}
\def \deg  {\mathop{\rm deg}\nolimits}

\def \Spec {\mathop{\rm Spec}\nolimits}

\def \socle {\mathop{\rm socle}}
\def \dim{\mathop{\rm dim}\nolimits}
\def \codim{\mathop{\rm codim}\nolimits}

\def \rank {\mathop{\rm rank}\nolimits}
\def \iso {\cong}
\def \tensor {\otimes}
\def \Sym {\mathop{\rm Sym}\nolimits}

\def \Pic {\mathop{\rm Pic}}
\def \Hom {\mathop{\rm Hom}\nolimits}
\def \Tor {\mathop{\rm Tor}\nolimits}
\def \Ext {\mathop{\rm Ext}\nolimits}
\def \Sec {\mathop{\rm Sec}\nolimits}

\def \th {{^{\rm th}}}
\def \st {{^{\rm st}}}
\def \red {{\rm red}}

\def \P {{\bf P}}
\def \H {{\rm H}}
\def \F {{\cal F}}
\def \K {{\cal K}}
\def \O {{\cal O}}
\def \I {{\cal I}}

\def \E {{\cal E}}
\def \h {{\rm h}}
\def \m {{\goth m}}

%%
% A small macro producing a box: #1 is the size of the
% box, #2 is what's in it !
%%

%

\def \tvi {\vrule height 12pt depth 7pt width 0pt}
\def \cc#1{\hfill\kern .7em#1\kern .7em\hfill}
%------------------------------------------
% end of begin.tex
%------------------------------------------

\input diagrams.tex

%
% Forward references
%
%-----------------------------------------
\forward {history}{Section}{1}
\forward {definitions}{Section}{2}
\forward {Gale-adjoint}{Section}{5}
\forward {examples}{Section}{3}
\forward {determinantal}{Section}{6}
\forward {gorenstein} {Section}{7}
\forward {linalg}{Section}{8}
\forward {classification} {Section}{9}
\forward {r+3} {Section} {4}
\forward {gale-scheme} {Definition} {2.1}
\forward {Goppa}{Corollary}{3.2}
\forward {stability}{Proposition}{8.10}
\forward {bath-orth-bases}{Theorem}{8.1}
\forward {determinants-and-Veronese}{Theorem}{6.1}
\forward {locally-Gorenstein}{Proposition}{6.2}
\forward {vero-gale}{Theorem}{7.2}
\forward {fail}{Theorem}{7.3}
\forward {quadric-fail}{Corollary}{7.4}
\forward {degen-self-assoc}{Proposition}{7.6}
\forward {clebsch}{Example}{5.12}
%-----------------------------------------

%\showlabels
%\rightline{\timestamp}
\bigskip

\centerline{\hd The Projective Geometry of the Gale Transform}
\medskip
\centerline {by}
\smallskip
\centerline {\bf David Eisenbud and Sorin Popescu
\footnote{$^{*}$}{\rm Both authors
are grateful to the NSF for
support during the preparation of this work.
\hfill\break\indent
It was in joint work with David Buchsbaum
that the first author first became familiar with the notion
of a Gorenstein ring. A large part of this 
paper (``self-associated sets'') is concerned, from an algebraic point
of view, with the classification and study of a special type of
Gorenstein ring, generalizing some of the examples found
in that joint work. It is with especial pleasure that we dedicate
this paper to David.}
}

\bigskip\bigskip
{\narrower
\abst The Gale transform, an involution on sets of points in 
projective space, appears in a multitude of guises,
in subjects as diverse as optimization, coding theory,
theta-functions, and recently in our proof that certain general sets
of points fail to satisfy the minimal free resolution conjecture (see
Eisenbud-Popescu [1996]).  In this paper we reexamine the Gale
transform in the light of modern algebraic geometry. We give a
more general definition, in the context of finite (locally) Gorenstein
subschemes. We put in modern form a number of the more remarkable
examples discovered in the past, and we add new constructions and
connections to other areas of algebraic geometry. We generalize Goppa's
theorem in coding theory and we give new applications to 
Castelnuovo theory. We give references to
classical and modern sources.\par}

\bigskip
\noindent{\bf Contents}
\medskip
\noindent I. The Gale Transform
\smallskip
\ref{history}: History

\ref{definitions}: The Scheme-theoretic Definition

\ref{examples}: A Generalized Goppa Theorem

\ref{r+3}: Castelnuovo's $r+3$ Theorem

\ref{Gale-adjoint}: The Gale Transform and Canonical Modules

\ref{determinantal}: The Gale Transform  of Determinantal Schemes
\medskip
\noindent II. Self-associated Schemes and Gorenstein Rings
\smallskip
\ref{gorenstein}: Gorenstein and Self-associated Schemes

\ref{linalg}: Linear Algebra and Self-association

\ref{classification}: Classification of Self-Associated
Schemes in Small Projective Spaces

\bigskip\goodbreak

Let $r,s$ be positive integers, and let $\gamma=r+s+2$.
The classical {\it Gale transform} is an involution that takes a
(reasonably general) set $\Gamma\subset \P^r$ of $\gamma$ 
labeled points in a projective space
$\P^r$ to a set $\Gamma'$ of
$\gamma$ labeled points in $\P^s$, defined up to a 
linear transformation of $\P^s$.
Perhaps the simplest (though least geometric)
of many equivalent definitions is this:
if we choose homogeneous coordinates so that 
the points of $\Gamma\subset\P^r$ have as coordinates the rows of the matrix
$$\pmatrix{ 
\tvi I_{r+1}\cr
\noalign{\hrule}  
\tvi A\cr},$$
where $I_{r+1}$ is an $(r+1)\times (r+1)$ identity matrix and $A$ is
a matrix of size $(s+1)\times (r+1)$, then the Gale transform
of $\Gamma$ is the set of points $\Gamma'$ whose homogeneous coordinates in 
$\P^s$ are the rows of the matrix 
$$\pmatrix{\tvi A^T\cr \noalign{\hrule}\tvi I_{s+1}\cr},$$
where $A^T$ is the transpose of $A$.

It is not obvious from this definition that the Gale transform
has any ``geometry'' in the classical projective sense.
Here are some examples that suggest it has:
\medskip

\noindent{$a)$ $r=1$:} The Gale transform of a set of $s+3$ points
in $\P^1$ is the corresponding set of $s+3$ points on the rational
normal curve that is the $s$-uple embedding of $\P^1$ in $\P^s$. Conversely,
the Gale transform of any $s+3$ points in linearly general position
in $\P^s$ is the same set in the $\P^1$ that is the unique 
rational normal curve through the original points. See \ref{Goppa}
and the examples following it.

\smallskip
\noindent{$b)$ $r=2, s=2$:} 
There are two main cases: a complete intersection of 
a conic and a cubic is its own Gale transform (a ``self-associated
set''). On the other hand, if $\Gamma$ consists of 6 points not on 
a conic, then the Gale transform of $\Gamma$ is the image of 
the 5 conics through 5 of the 6 points of $\Gamma$ via the Cremona
transform that blows up the 6 points and then blows down the conics.
See \ref{clebsch}.

\smallskip
\noindent{$c)$ $r=2, s=3$:} A set $\Gamma$ of 7 general points in $\P^3$ lies
on 3 quadrics, which intersect in 8 points. The Gale transform
of $\Gamma$ is the projection of $\Gamma$ from the eighth point.
Again, see the examples following \ref{Goppa}.

\medskip
\ref{history}  contains what we know of the 
history of the Gale transform, including work of Pascal,
Hesse, Castelnuovo, Coble, Dolgachev-Ortland, and Kapranov.

In \ref{definitions} we introduce a general definition of the Gale
Transform as an involution, induced by Serre duality, on the set of
linear series on a finite Gorenstein scheme over a field (here as
always in this paper, Gorenstein means locally Gorenstein). This
language turns out to be very convenient even in the classical
(reduced) case.  The main result of this section is an extension of
the Cayley-Bacharach Theorem for finite complete intersections (or, more
generally, finite arithmetically Gorenstein schemes) to finite locally
Gorenstein schemes. It interprets the failure of a set of points to
impose independent conditions on a linear series as a condition on the
Gale transform.

The next sections treat basic properties of the Gale transform,
and examples derived from them, such as the ones above.
\ref{examples} is devoted to an extension of a famous theorem of Goppa in
coding theory. A linear code is essentially a set of points in
projective space, and the dual code is its Gale transform.  Goppa's
theorem asserts that if a linear code comes from a set of points on a
smooth linearly normal curve, then the dual code lies on another image
of the same curve. Examples a) and c) above are special
cases. We show how to extend this theorem (and its
scheme-theoretic generalization) to sets of points contained in
certain other varieties, such as ruled varieties over a curve. Using
these results we exhibit some of the classical examples of the Gale
transform, and provide some new ones as well; for instance 9 general
points of $\P^3$ lie on a smooth quadric surface, which is a ruled
variety in two different ways. It follows from our theory that the
Gale transform, which will be 9 general points of $\P^4$, lies on two
different cubic ruled surfaces. In fact, we show that the 9 general
points are the complete intersection of these two surfaces.

In \ref{r+3} we use the Gale transform to give a simple proof of
Eisenbud-Harris' generalization to schemes of Castelnuovo's lemma that
$r+3$ points in linearly general position in $\P^r$ lie on a unique
rational normal curve in $\P^r$. We also prove a similar result on
when finite schemes in linearly general position lie on higher
dimensional rational normal scrolls, and when these scrolls may be
taken to be smooth. Our method provides a simple proof (in many cases) 
for a result of Cavaliere-Rossi-Valla [1995].

In \ref{Gale-adjoint} we show that if $\Gamma\subset\P^r$ and
$\Gamma'\subset\P^s$ are related by the Gale transform, then
the canonical modules $\omega_\Gamma$ and $\omega_{\Gamma'}$
are related in a simple way. This is the  idea exploited
in Eisenbud-Popescu [1996] to study the minimal free
resolutions associated to general point sets $\Gamma\subset\P^r$
and in particular to disprove the Minimal Resolution Conjecture.
As an application we exhibit an example due to Coble, connecting
the Gale transform of 6 points in the plane with the Clebsch transform
(blow up the six points, blow down the proper transforms of the
conics through five of the six.)

One family of examples that does not seem to have been considered
before are the determinantal sets of points. In \ref{determinantal} we
describe a novel relationship, expressed in terms of the Gale
transform of Veronese re-embeddings, between the zero-dimensional
determinantal varieties defined by certain ``adjoint'' pairs of
matrices of linear forms. 

A major preoccupation of the early work on the Gale transform was the
study of ``self-associated'' sets of points, that is, sets of points
$\Gamma\subset\P^r$ that are equal to their own Gale transforms (up to
projective equivalence, of course).  This notion only
applies to sets of $2r+2$ points in $\P^r$, since otherwise the Gale
transform doesn't even lie in the same space.  For example, 6 points
in the plane are self-associated iff they are complete intersection of
a conic and a cubic (and this is the essential content of Pascal's
``Mystic Hexagram'').  It turns out that this is indeed a natural
notion: under mild non-degeneracy assumptions $\Gamma$ is
self-associated iff its homogeneous coordinate ring is Gorenstein!
\ref{gorenstein} is devoted to a study of 
self-association and a generalization: Again under mild extra
hypotheses, a Gorenstein scheme $\Gamma\subset\P^r$ has
Gorenstein homogeneous coordinate ring iff the Gale transform
of $\Gamma$ is equal to a Veronese transform of $\Gamma$.
We review the known geometric constructions of self-associated
sets, and add a few new ones. It would be interesting to 
know whether the list contains any families of Gorenstein
ideals not yet investigated by the algebraists.

In \ref{linalg} we continue the study of self-associated sets of
points, showing how they are related to nonsingular bilinear forms on
the underlying vector space of $\P^r$.  A classical result states that
self-associated sets correspond to pairs of orthogonal bases of such a
form.  We say what it means for a non-reduced scheme to be the ``union
of two orthogonal bases'', and generalize the result correspondingly.
We also reprove and generalize some of the other classic results on
self-associated sets, showing for example that the variety
of self-associated sets of labeled points in $\P^r$ is isomorphic to
an open set of the variety of complete flags in $\P^r$, a result of
Coble and Dolgachev-Ortland (see the references below).

It is interesting to ask, given a set of $\gamma$ points in $\P^r$ with 
$\gamma<2r+2$, whether it can be extended to a set of points
of degree $2r+2$ with Gorenstein homogeneous coordinate ring---indeed,
such questions arise implicitly in our work on free resolutions
[1996]. From the theory developed in \ref{linalg} we are able to
give interesting information in some cases. For example, we show
that a set of 11 general points in $\P^6$ can be completed to 
a set of $14$ points with arithmetically Gorenstein
homogeneous coordinate ring. We show that (although the extension
is not unique) the three points added span a plane that is uniquely
determined. This plane appears as the ``obstruction'' to the truth
of the minimal resolution conjecture for 11 points in $\P^6$, as
treated in our paper [1996].

In \ref{classification} we continue with self-associated sets,
and  describe what is known about the  classification of
small dimensional projective spaces, up to $\P^5$.

We thank Joe Harris and Bernd Sturmfels for introducing us to the Gale
transform and to Lou Billera, Karen Chandler and Tony Geramita for
useful discussions.

\par

\section{history} History

Perhaps the first result that belongs to the development of the
Gale transform is the theorem of Pascal (from his ``Essay Pour Les
Coniques'' from 1640, reproduced in Struik [1969])
that the vertices of two triangles circumscribed around the same conic
lie on another conic. As we shall see in \ref{linalg},
this is a typical result
about sets of points that are Gale transforms of themselves
(``self-associated
sets'').
Hesse, in his Dissertation and Habilitationschrift in K\"onigsberg [1840]
(see the paper in Crelle's Journal [1840] and the reprinted Dissertation in
Werke [1897]), found an analogue of Pascal's result (see also
Zeuthen [1889, p. 363]) that held for 8 points in
three dimensional space and gave various applications.
Some of Hesse's results were made clearer and
also extended by von Staudt [1860], Weddle [1850],
Zeuthen [1889], and Dobriner [1889] (see also p.152,
Tome 3, of the Encyclopedie [1992])

The step to defining the Gale transform itself in the
corresponding cases of 6 points in $\P^2$ was taken by Sturm [1877],
and extended by Rosanes [1880, 1881].
More important is the realization by Castelnuovo in [1889]
that one could do the same sort of thing for $2r+2$ points
in $\P^r$ in general. He called two  sets of $2r+2$ points that are
Gale transforms of one another ``gruppi associati di punti''.
Castelnuovo, who refers to Sturm and Rosanes but seems unaware
of Hesse's work, gives the following geometric definition:
\par\medskip

{\sl Two sets $\Gamma$ and $\Gamma^\prime$, each of $2r+2$ labelled points
in  $\P^r$ and ${(\P^r)}^\ast$ respectively, are defined to be
associated when there exist two
simplices $\Delta$ and $\Delta^\prime$ in $\P^r $ such that
the points of $\Gamma$ are projective with the $2r+2$ vertices
of $\Delta$ and $\Delta^\prime$, while the points of $\Gamma^\prime$ are
projective with the  $2r+2$ facets of the two simplices (each facet
being labeled by the opposite vertex)}.
\par\medskip

Castelnuovo was primarily interested in the case when a set of
points is self-associated;
this is the case of Pascal's 6 points on a conic, for example.
As we shall see, a general set of $2r+2$ points is associated to itself
iff its homogeneous coordinate ring is Gorenstein. It
is interesting that the stream of
work that lead Castelnuovo in this direction has the same source in
Pascal's theorem as the stream that lead to the Cayley-Bacharach Theorem and
its ramifications, another early manifestation of the Gorenstein
property (see Eisenbud-Green-Harris [1996] for a discussion).

The first one to have studied the Gale transform of a set of $\gamma$
points in $\P^r$ without assuming $\gamma=2r+2$ seems to have been
Coble, and we begin with his definition. As before, we may represent
an ordered set of $\gamma$ points $\Gamma\subset\P^r$ 
by a $\gamma\times (r+1)$ matrix of homogeneous
coordinates, though this involves some choices. (To make the symmetry
of the relation of $\Gamma$ and its Gale transform better visible, we
will no longer insist, as above, that the first part of the matrix is
the identity.)

\definition{gale-naive}  Let $k$ be a field, and let $r,s\geq 1$
be integers.  Set $\gamma=r+s+2$, and let
$\Gamma\subset\P^r,\, \Gamma'\subset\P^s$
be ordered nondegenerate sets of $\gamma$
points represented by $\gamma\times (r+1)$
and $\gamma\times(s+1)$ matrices $G$ and $G'$, respectively.
We say that $\Gamma'$ is the {\it Gale transform}
of $\Gamma$ if there exists a nonsingular diagonal
$\gamma\times\gamma$ matrix $D$
such that $G^{\rm T}\cdot D\cdot G'=0$

Put more simply, the Gale transform of a set of points represented
by a matrix of homogeneous coordinates $G$ is the set of points
represented by the kernel of $G^{\rm T}$ (the diagonal matrix
above is necessary to avoid the dependence on the choices of homogeneous
coordinates). Note that the Gale transform is really only defined up
to automorphisms of the projective space. Since we are going
to give a more general modern definition in the next section,
we will not pause to analyze this version further.

This definition is related to the one given in the 
introduction by the identity
$$
\pmatrix{A& \mid & I_{s+1}}
\pmatrix{-I_{r+1}& 0 \cr 0 &I_{s+1}}
\pmatrix{\tvi I_{r+1}\cr\noalign{\hrule} \tvi A}
=0.
$$

In a remarkable series of papers ([1915, 1916, 1917, 1922]) in the
early part of this century Coble gave the definition above, gave
applications to
theta functions and Jacobians of curves, and described many amazing
examples. Although Coble used the same term for
the ``associated sets'' as Castelnuovo, he doesn't mention Castelnuovo
or any of the other references given above, leaving us to wonder how
exactly he came to the idea.  (The related paper of Conner [1911],
often quoted by Coble, doesn't mention Castelnuovo either.)

Castelnuovo's work on self-associated sets of points, on the
other hand, was continued by Bath [1938], Ramamurti [1942]
and Babbage [1948], but
they seem ignorant of Coble; perhaps the old and new worlds
were too far apart.

Similar ideas in the affine case were developed, apparently
without any knowledge of this earlier work, by Whitney [1940]
and Gale [1963], the latter in the study of polytopes.
As a duality theory for polytopes, and in linear programming,
it has had a multitude of applications.
The names ``Gale Transform'' and ``Gale diagram'' are
well-established in these fields, and in the absence of
a more descriptive term than ``associated points'' in algebraic
geometry we have adopted them.

Another important group of applications was initiated by Goppa [1970,1984].
In coding theory the Gale transform is the passage from a code to
its dual, and Goppa proved that a code defined by a set of rational
points on an algebraic curve was dual to another such code. See
\ref{Goppa} for a generalization.

Dolgachev and Ortland [1988] give a modern exposition of the geometric
theory of the Gale transform. These authors treat many topics covered
by Coble.  Their main new contribution is the use of geometric
invariant theory to extend the definition of the Gale transform to a
partial compactification of the ``configuration space'' of sets of
general points. In a similar vein, Kapranov [1993]
shows that the Gale transform extends from general sets of 
points to the Chow compactification. (A similar result can
be deduced from our description of the Gale transform as an
operation on linear series, since the Chow compactification
is an image of the set of linear series on a reduced set of
points containing a given generator of the line bundle.)

Although we will not pursue it here, there is a possible
extension of the theory, suggested to us by Rahul Pandharipande,
which deserves mention.
One may easily extend \ref{gale-naive} to ordered collections
of linear subspaces:  the Gale transform of a collection
of $\gamma$ linear spaces of dimension $d$ in $\P^r$ will be
a collection of $\gamma$  linear spaces of dimension also
$d$ in $\P^{\gamma(d+1)-r-2}$. (For each space, choose an independent
set of points spanning it. Take the Gale transform 
of the union of these collections of points. The spaces
spanned by the subsets corresponding to the original spaces
make up the Gale transform.) Thus the Gale transform of a set
of 4 lines in $\P^3$ would be again a set of 4 lines in $\P^3$.
This definition can be shown to be independent of the choice of
frames  and generalized to finite Gorenstein subschemes
of Grassmannians in the spirit of \ref{gale-scheme} bellow. 
It might be interesting to understand its geometric significance,
at least in such simple examples as that of the four lines
above.

\section{definitions} The Scheme-theoretic definition

Throughout this paper $\Gamma$ will denote
a Gorenstein scheme, finite over a field $k$ (we will usually
just say a {\sl finite Gorenstein scheme\/.}) We recall that
a subscheme of $\P^n$ is {\it nondegenerate\/} if it is 
not contained in any hyperplane.

The Gale transform is the involution
on the set of linear series on $\Gamma$ induced
by Serre duality. In more detail:

Let  $\Gamma$ be  a finite Gorenstein scheme.
let $L$ be  a line bundle on $\Gamma$, and consider
the canonical ``trace'' map $\tau:\ \H^0(K_{\Gamma})\rTo^\tau k$
provided by Serre duality.  The composition
of $\tau$ with the multiplication map
$$
\H^0(L)\tensor_k \H^0(K_{\Gamma}\otimes L^{-1}) \rTo \H^0(K_{\Gamma}) \rTo k,
$$
is a perfect pairing  between
$\H^0(L)$ and $\H^0(K_{\Gamma}\otimes L^{-1})$.  For any subspace
$V\subset \H^0(L)$  we write
$V^\perp \subset \H^0(K_{\Gamma}\otimes L^{-1})$ for the annihilator
with respect to this pairing.

\definition{gale-scheme} Let $\Gamma$ be  a finite Gorenstein scheme,
let $L$ be a line bundle on $\Gamma$, and let
$V\subset \H^0(L)$ be a subspace. The Gale Transform of the
linear series $(V,L)$ is the linear series
$(V^\perp,K_{\Gamma}\otimes L^{-1})$.

In Coble's work [1915, 1922] the following observation,
always in the reduced case and with $\deg \Gamma_1=r+1$,
is used as the foundation of the theory:

\proposition {basic-principle} Let $\Gamma$ be
a finite Gorenstein scheme of degree
$r+s+2$.
Let $\Gamma_1\subset \Gamma$ be a
subscheme, and let $\Gamma_2$ be the residual scheme to
$\Gamma_1$.  Let $(V,L)$ be a linear series of (projective)
dimension $r$ on $\Gamma$, and suppose that $W\subset
\H^0(K_\Gamma\tensor L^{-1})$ is its Gale transform.  The
failure of $\Gamma_1$ to span $\P(V)$ (that is, the codimension
of the linear span of the image of $\Gamma_1$ in $\P(V)$)
is equal to the failure of $\Gamma_2$ to impose independent
conditions on $W$.

\proof Consider the diagram with short exact row and 
column
$$\diagram[small,midshaft]
&&\H^0(L_{|\Gamma_1})&\\
&\ruTo^a&\uTo&\\
V&\rTo&\H^0(L)&\rTo&W^*\\
&&\uTo&\ruTo_b\\
&&\H^0(\I_{\Gamma_1}L).
\enddiagram$$
A diagram chase shows that $\ker(a)\iso \ker(b)$.  But the failure of
$\Gamma_1$ to span $\P(V)$ is the dimension of the kernel of $a$, and
since $\H^0(\I_{\Gamma_1}L)^*=\H^0(L\mid_{\Gamma_2})$, the failure of
$\Gamma_2$ to impose independent conditions on $W$ is the dimension of
the cokernel of $b^*$.\Box

\problem{} Coble often uses this result as follows: he gives
some transformation taking
a general set of $r+s+2$ labeled points in $\P^r$ to a general set of
$r+s+2$ labeled points in $\P^s$, definable by rational functions. He
then proves that it takes a set of points whose first $r+1$ elements are
dependent to a set of points whose last $s+1$ elements are
dependent. He then claims that in consequence this transformation must
agree with the Gale transformation. Coble establishes many of the
examples described below in this way.  Can this argument be made
rigorous?
\medskip

As a corollary of \ref{basic-principle} we deduce a characterization
of base-point-free and very ample
linear series in terms of their Gale transforms. It will
usually be applied to the Gale transform of a known linear
series $V$, so we formulate it for $W=V^\perp$.

\corollary{bpf}  Let $\Gamma$ be a finite
Gorenstein scheme over a field $k$ with algebraic closure 
$\overline{k}$, let $L$ be a line bundle on $\Gamma$,
and let $V\subset\H^0(L)$ be a linear series. 
Let $W=V^\perp\subset\H^0(K_\Gamma\tensor L^{-1})$ 
be the Gale transform of $(V,L)$.
\item{a)} The series $W$ is base point free iff 
no element of $\overline k\tensor V$ vanishes on a codegree 1 subscheme of 
$\overline k\tensor \Gamma$. 
\item{b)} The series $W$ is very ample iff no element of
$\overline{k}\tensor V$ vanishes on a codegree 2 subscheme
of $\overline{k}\tensor\Gamma$.

For example, suppose that $\Gamma$ is reduced  over $k=\overline k$ and 
$\Gamma$ is embedded by $V$ into $\P(V)$. The series $W$
is base point free iff no hyperplane in $\P(V)$ contains all but one
point of $\Gamma$; and $W$ is very ample iff no hyperplane in $\P(V)$
contains all but two points of $\Gamma$.

In the case of basepoint freeness, a criterion can be given which
does not invoke $\overline k$. The reader may check that condition
$a)$ is equivalent to the statement that no proper subscheme $\Gamma'$
of degree $\gamma'$ in $\Gamma$ imposes only $\gamma'-\dim(W)$
conditions, the smallest possible number, on $W^\perp$.

\proof Both statements reduce immediately to the
case $k=\overline k$. Then $W$ is basepoint free iff every
degree 1 subscheme of $\Gamma$ imposes one condition on $W$. By
\ref{basic-principle}, this occurs iff every codegree 1 subscheme
of $\Gamma$ spans $\P(V)$, that is, iff no element of $V$ 
vanishes on a codegree 1 subscheme.  Similarly, $W$ is 
very ample iff every degree 2 subscheme imposes 2 conditions
on $W$, and again the result follows from \ref{basic-principle}.
\Box

\smallskip
\ref{basic-principle} may be seen as a generalization
of the Cayley-Bacharach theorem from the case of arithmetically
Gorenstein schemes to the case of (locally) Gorenstein schemes.  (See
Eisenbud-Green-Harris [1996] for historical remarks on this other forerunner
of the Gorenstein notion.)  To exhibit this aspect we first explain 
how to find the Gale transforms of series cut out by hypersurfaces of
given degree.

Recall that if $\Gamma\subset\P(U)$ is any finite scheme, with
homogeneous coordinate ring $S_\Gamma$, then
the {\it canonical module\/} of $\Gamma$ is the $S_\Gamma$-module
$$
\omega_\Gamma = \Ext_S^{\dim(U)-1}(S_\Gamma, S(-\dim(U))),
$$
where $S$ denotes the symmetric algebra of $U$,
the homogeneous coordinate ring of $\P(U)$. 
The sheaf on $\Gamma$ associated to $\omega_\Gamma$ is the 
dualizing sheaf $K_\Gamma$.

\proposition{hypersurface-series} Let $\Gamma\subset\P(U)$ be a finite
Gorenstein scheme, and let $U^d$ be the series cut out on $\Gamma$
by the hypersurfaces of degree $d$ in $\P(U)$. The Gale transform
of $U^d\subset\H^0(\O_\Gamma(d))$ is the image of ${(\omega_\Gamma)}_{-d}$ in 
$\H^0(K_\Gamma(-d))$.

\proof From the exact sequence 
$0\rTo\I_\Gamma \rTo \O_{\P(U)}\rTo\O_\Gamma\rTo 0$ 
we get the sequence
$$
\dots\rTo
\H^0(\O_{\P(U)}(d))\rTo^\alpha \H^0(\O_\Gamma(d))\rTo \H^1(\I_\Gamma(d))\rTo 0
$$
which identifies $(U^d)^\perp$ with $\H^1(\I_\Gamma(d))^*$.
By duality, $\omega_\Gamma$ is the module
$$\bigoplus_{e=-\infty}^\infty \H^1(\I_\Gamma(-e))^*,$$ with degree $-d$
part equal to $\H^1(\I_\Gamma(d))^*$, so we are done.\Box

For any finite scheme $\Gamma\subset\P(V)$, let $a(\Gamma)$ be the
largest integer $a$ such that $\Gamma$ fails to impose independent
conditions on forms of degree $a$. For example, if $\Gamma$ is a
complete intersection of forms of degrees $d_1,\ldots,d_c$, then
$a(\Gamma)= \sum (d_i-1)-1$. The scheme $\Gamma$ is {\it
arithmetically Gorenstein\/}---that is, the homogeneous coordinate
ring $S_\Gamma$ is Gorenstein---iff there is an isomorphism
$S_\Gamma(a)\iso \omega_\Gamma$ for some $a$, and it is easy to see
that then $a=a(\Gamma)$. See for example Bruns-Herzog [1993],
Proposition 3.6.11 and Definition 3.6.13.

\corollary{gale-for-gor}  Suppose that the finite scheme
$\Gamma\subset\P(U)$ is arithmetically Gorenstein.  For any integer
$d\leq a(\Gamma)$, the linear series of forms of degree $d$ on
$\Gamma$ is the Gale transform of the linear series of forms of degree
$a(\Gamma)-d$.

\proof Use \ref{hypersurface-series} and the fact that
$\omega_\Gamma = S_\Gamma(a)$.\Box

The Cayley-Bacharach theorem is now the special case of 
\ref{basic-principle} in which $\Gamma$ is arithmetically
Gorenstein in some embedding in a $\P^n$,
and the linear series involved is induced by forms of  small
degree on this projective space.

\corollary{Cayley-Bacharach} {\bf (Cayley-Bacharach for arithmetically
Gorenstein schemes)} Suppose that
$\Gamma\subset \P(U)$ is a finite arithmetically Gorenstein 
scheme, and let $\Gamma_1$, $\Gamma_2$ be mutually residual
subschemes of $\Gamma$. For any integer $d< a(\Gamma)$,
the failure of $\Gamma_1$ to impose independent
conditions on forms of degree $d$ is
equal to the number of forms of degree
$a(\Gamma)-d$ vanishing on $\Gamma_2$.

\proof Apply \ref{basic-principle} and \ref{gale-for-gor}.
\Box

This includes the classic version:

\corollary{Chasles} ({\bf Chasles}) If a set $\Gamma_1$ of
8 points in $\P^2$ lies in the complete
intersection $\Gamma$ of two cubics, then any cubic 
vanishing on $\Gamma_1$
vanishes on $\Gamma$.

\proof  In this case $a(\Gamma)=3$. 
Since the number of forms of degree $0$ vanishing
on the empty set (respectively any one-point set) is 1 (respectively
0), the nine points of $\Gamma$ impose dependent conditions on cubics,
while any 8-point subset imposes independent conditions on cubics.\Box

\section{examples} A Generalized Goppa Theorem

This section is devoted to an extension of Goppa's classical result
on the duality of algebro-geometric codes.

\theorem{vbdle-gale}
Let $\Gamma$ be a zero-dimensional  Gorenstein scheme with a finite
map to a locally Gorenstein base scheme $B$ of dimension $c$, and let
${\cal O}_{\Gamma}(1)$ be a line bundle on $\Gamma$.  Suppose
that
$$
0\rTo \E_{c}\rTo\E_{c-1}\rTo \ldots \rTo \E_{0}\rTo\O_\Gamma(1)\rTo 0
$$
is a resolution of $\O_\Gamma(1)$ by locally free sheaves on $B$,
and hence that
$$
{K}_{\Gamma}(-1) \ = \
\coker \bigl[{\cal H}om(\E_{c-1},K_{B})
\rTo\relax {\cal H}om(\E_{c},K_{B})\bigr].
$$
If
$$
\H^{i+1}(\E_i)=\H^{i+1}(\E_{i+1})=0,
\quad {\sl for\ all}\quad 0\le i\le c-2,
$$
then the induced sequence
$$\eqalign{
\H^{0}(\E_{1})\rTo \H^{0}(\E_{0})\rTo
&\bigl[\H^{0}(\O_{\Gamma}(1))=
\H^{0}({K}_{\Gamma}(-1))^\ast\bigr] \rTo   \cr
&\H^{0}({\cal H}om(\E_{c},K_{B}))^\ast
\rTo \H^{0}({\cal H}om(\E_{c-1},K_{B}))^\ast
}$$
is exact. In particular, if $\Gamma$ is embedded in
$\P^{r}$ by the linear
series which is the image of $\H^{0}(\E_{0}) \rTo \H^{0}({\O}_{\Gamma}(1))$,
so that $\Gamma$ lies on the image of $\P(\E_{0})$ embedded by the
complete linear series $|{\O}_{\P(\E_{0})}(1)|$, then the
Gale transform of $\Gamma$ is defined by the image of
$\H^{0}(\E_{c}^*\otimes K_B)\rTo \H^{0}({K}_{\Gamma}(-1))$, and  lies on the
image of $\P(\E_{c}^*\otimes K_B)$ mapped by the
complete linear series $|{\O}_{\P(\E_{c}^*\otimes K_B)}(1)|$.

\proof  Break up the given resolution
into short exact sequences
$$
0\rTo \K_{i+1}\rTo \E_i\rTo \K_{i}\rTo 0,\quad 0\le i\le c-1,
$$
and then use the vanishings in the hypothesis to obtain
$$\eqalign{
\coker (\H^{0}(\E_{0})\rTo \H^{0}({\O}_{\Gamma}(1)))=&\H^1(\K_1) =
\ldots \cr
&= \H^{c-1}(\K_{c-1})= \ker (\H^c(\K_c) \rTo \H^c(\E_{c-1})).
}$$
Identifying $\K_c$ with $\E_c$ and using Serre duality,
we get the asserted result. \Box

\smallskip
 In the special case when $\Gamma\subset B=\P^r$,
\ref{vbdle-gale} is the special case $d=1$ of
\ref{hypersurface-series}.

\corollary{Goppa} (Goppa Duality) Let $B$ be a locally
Gorenstein curve, embedded in
$\P^{r}$ by the complete linear series associated to a line bundle
${\O}_B(H)$, and let $\Gamma\subset B\subset \P^{r}$ be a Cartier
divisor on the curve $B$.  The Gale transform of
$\Gamma$ lies on the image of $B$ under the complete
linear series associated to ${\O}_B(K_B-H+\Gamma)$.

\proof Apply \ref{vbdle-gale} to the resolution
$$
0\rTo \O_B(H-\Gamma)\rTo\O_B(H)%
\rTo\O_{\Gamma}(1)\rTo 0.
$$
\Box

This result is essentially  due to Goppa [1970], [1984], and
expresses the duality among the algebro-geometric
codes bearing his name (see e.g. van Lint-van der Geer
[1988] for more details). From  \ref{Goppa}
we may immediately derive the following consequences:
\item{$\bullet$} If a curvilinear finite scheme of degree $\gamma=r+s+2$
lies on a rational normal curve in $\P^r$, then its Gale transform lies also
on a rational normal curve in $\P^s$. The analogous statement
also holds for finite subschemes of an elliptic
normal curve. See Coble [1922] for the statement in the reduced
case.
\item{$\bullet$} A set of $\gamma=2g-2$ points in $\P^{g-2}$
which is the hyperplane section of a canonical curve of genus
$g$ is its own Gale transform.
\item{$\bullet$} Let $\Gamma$ be a set of seven general points in $\P^2$,
and let $E\subset\P^2$ be a smooth plane cubic curve passing through
$\Gamma$.  Write $h$ for the hyperplane divisor of $E\subset \P^2$.
By \ref{Goppa}, the Gale
transform $\Gamma^\prime$  of $\Gamma$ is the image of $\Gamma$
via the re-embedding of $E$ as an elliptic normal
quartic curve $E^\prime\subset\P^3$ with hyperplane divisor
$H=\Gamma-h$. Similarly, $\Gamma$ is obtained from 
$\Gamma'$ by the linear series $\Gamma'-H$. If we write 
$\Gamma' = 2H-p$, so that the 3 quadrics containing $\Gamma'$
intersect in $\Gamma'+p$, then we see that $\Gamma'-H=H-p$, so
that $\Gamma$ is obtained from $\Gamma'$ by projection from $p$.

\corollary{scroll-gale}  Let $B$ be a locally Gorenstein
curve, and let $\E$ be a vector bundle over $B$.
Let $\Gamma$ be a zero dimensional Gorenstein
subscheme of the ruled variety $X = \P(\E)$, and
assume that $\Gamma$ is embedded
in $\P^{r}$ by the restriction of the complete series
$|{\O}_{\P(\E)}(1)|$. Assume
that $\O_{\P(\E)}(1){\mid_\Gamma}$ is very ample. 
Then the Gale transform of $\Gamma$ lies on the image of the
ruled variety $X' = \P({(\E^\prime)}^\ast\otimes K_{B})$, mapped
by $|\O_{\P({(\E^\prime)}^\ast\otimes K_{B})}(1)|$, where
the vector bundle $\E^\prime$ is defined as the kernel of
the natural epimorphism
$$
0\rTo\E^\prime\rTo\E\rTo\O_{\Gamma}(1)\rTo 0.
$$
If $X$ is a ruled surface, that is, $\rank(\E) = 2$, then
${(\E^\prime)}^\ast\otimes K_{B}\cong
\E'\otimes\det(\E')^\ast\otimes K_{B}$,
and hence $X^\prime$ is  the elementary transform
of $X$ along the scheme $\Gamma$.

\proof Apply \ref{vbdle-gale} to the resolution
$$
0\rTo\E^\prime\rTo\E\rTo\O_{\Gamma}(1)\rTo 0,
$$
\Box

\example{nine3} A smooth quadric surface $Q$
in $\P^3$ can  be regarded in two ways as a ruled surface over
$\P^1$, hence we deduce that 9 general points in $\P^4$ lie in the
intersection of two rational cubic scrolls in $\P^4$. 
The 9 points are actually the complete intersection of the two
scrolls. 

To see this, let $\Gamma\subset Q\subset\P^3$ be a set of
9 general points. The ideal $I_\Gamma$ is 3-regular
so by Bertini the general cubic through
$\Gamma$ cuts out on $Q$ 
a general canonically embedded smooth curve of genus 4.
Such a curve has exactly two $g^1_3$s, namely
those cut out by the two rulings of the smooth quadric $Q$.
By \ref{Goppa}, the Gale transform of $\Gamma$ is a hyperplane
section of the re-embedding of $C$ in $\P^5$ via the linear
system $|\Gamma|$. In this embedding, each
$g^1_3$ of $C$ sweeps out a rational cubic threefold 
scroll $X_i$, $i=1,2$ (isomorphic to the Segre embedding of 
$\P^1\times \P^2$ into $\P^5$).  Each $X_i$ is determinantal, cut out
be the $2\times 2$ minors of a $2\times 3$ matrix with linear
entries. Since $C\subset\P^5$ is a general curve of degree
9 and genus 4 in $\P^5$, and since the complete intersection
of two general cubic scrolls is of this type,
$C=X_1\cap X_2$, and $\Gamma'$ is correspondingly 
the complete intersection of the two scrolls in $\P^4$.

If the points $\Gamma$ lie on a complete intersection of
a cubic with a singular quadric, then the two scrolls in $\P^4$
coincide, so $\Gamma'$ is not a complete intersection as above.
This leads us to the following formulation.

\problem{two-scroll-intersection} 
Suppose that $\Gamma\subset \P^3$ is a finite
Gorenstein scheme of degree  9 that lies on a unique quadric
surface, and suppose that this quadric is smooth.
Is the Gale transform of $\Gamma$ the
intersection of the corresponding pair of 
rational cubic scrolls in $\P^4$?

More generally the Gale transform of $5+r$ general points
on a smooth quadric surface in $\P^3$ is contained in the intersection
of two rational normal scrolls in $\P^r$.
\smallskip

We refer the reader to \ref{r+3} for other applications of
\ref{Goppa} and \ref{scroll-gale}.

\section{r+3} Castelnuovo's $r+3$ Theorem

A scheme-theoretic version of Castelnuovo's Lemma for $r+3$ points
in linearly general position (over an algebraically closed field)
was proved by Eisenbud and Harris [1992]:
\smallskip
A finite subscheme $\Gamma\subset\P^r_k$ of degree $r+3$
in linearly general position over an algebraically closed field $k$
lies on a unique rational normal curve.

\smallskip
As an  application of the Gale transform we give here a simpler
direct proof of this result.

It is well-known that any finite subscheme $\Gamma$ of a rational normal
curve over an arbitrary field $k$ 
is in {\it linearly general position\/}, in the sense that
any subscheme of $\Gamma$ that lies on a $d$-dimensional linear subspace
has degree $\leq d+1$. It is also {\it curvilinear\/}  (each local ring
$\O_{\Gamma,p}$ is
isomorphic to $F[x]/(x^n)$ for some $n$, and some field extension
$F$ of $k$) and {\it unramified\/}
(if $V$ is the vector space of linear forms on
$\P^r$, then the natural map $V\rTo \O_{\Gamma,p}$ has image the complement of
an ideal).  Thus the following
result characterizes subschemes of a rational normal curve:

\theorem{Castel-Gor}  Let $\Gamma$ be a finite scheme
geometrically in linearly general position in $\P^r_k$.
\item{a)} If $\deg \Gamma = r+3$, then $\Gamma$
lies on a rational normal curve iff
$\Gamma$ is Gorenstein.
\item{b)} If $\deg \Gamma\geq r+3$ and $\Gamma$ is  Gorenstein,
then $\Gamma$ is curvilinear and unramified.

\proof {\sl a)} If $\Gamma$ lies on a rational normal curve, then $\Gamma$ is
curvilinear, thus Gorenstein.
Conversely, suppose that $\Gamma$ is
Gorenstein and geometrically in linearly general position.  
By \ref{bpf}, the Gale transform of 
$\Gamma\subset \P^r_k$ is an embedding of $\Gamma$
as a subscheme $\Gamma'\subset \P^1_k$. By \ref{Goppa}, $\Gamma$
the Gale transform of $\Gamma'$, lies on a rational normal curve.
(One could also use \ref{vero-gale}: since any subscheme of
$\P^1_k$ is arithmetically Gorenstein, its Gale transform is
equal to its Veronese transform.)

{\sl b)} By part {\sl a)}, any subscheme of degree $r+3$ of $\Gamma$ is
curvilinear and unramified.  
It follows that every component of $\Gamma$ is too, and this implies
the desired result.\Box

\remark{} The condition of being geometrically
in linearly general position cannot be replaced by the condition of
being in linearly general position.  Let $F$ be a field of characteristic
$p$, and let $k=F(s,t)$.  Consider the Gorenstein
scheme $\Gamma= \Spec k(s^{1/p},t^{1/p})$.
Set $r=2p-3$,
and let $V\subset \O_\Gamma = k(s^{1/p},t^{1/p})$ be a $k$-subspace of
dimension
$r+1=2p-2$. If $p\geq 3$, then $V$ has dimension more than half the
dimension of $\O_\Gamma$, and thus $V$ is very ample. As
$\Gamma\subset\P^r_k=\P(V)$
has no proper subschemes at all, it is in linearly general position. But it
does NOT lie on any rational normal curve, since after tensoring with
$\overline k$
the local ring $\overline {k}\tensor \O_\Gamma \iso
\overline{k}[x,y]/(x^p,y^p)$
is not generated by one element over $\overline k$.

\medskip

To connect this result with the result for algebraically closed
ground fields proved by Eisenbud and Harris we use:

\theorem{castel-r+3} Let $\Gamma$ be a finite scheme
in linearly general position in $\P^r_k$.
If $k$ is algebraically closed and $\deg\Gamma\geq r+2$,
then $\Gamma$ is Gorenstein.

\example{}  The following shows that the
condition of algebraic closure cannot be dropped in
\ref{castel-r+3}.  Let $\Gamma$ be the scheme in $\P^2$
defined by the $2\times 2$ minors of the matrix
$$\pmatrix {
x&y&0\cr
-y&x&x^2
}.$$
$\Gamma$ is a finite scheme of degree 5, concentrated at the point $p$
defined by $x=y=0$. It is in linearly general position over
$\bf R$ (but not after base change to $\bf C$.)
It is not Gorenstein since
the matrix above gives a minimal set of syzygies locally at $p$.
In particular, it does not lie on a rational normal curve.
\medskip

\noindent{\sl Proof\ of \ref{castel-r+3}}.
The result amounts to a very special
case of Theorem 1.2 of Eisenbud-Harris [1992].  Here is a greatly simplified
version of the proof given there. See the original for further
classification, examples, and remarks.

Let $V\rTo \H^0(\O_\Gamma(1))$ be the map defining
the embedding of $\O_\Gamma$ in $\P^r_k=\P(V)$. By choosing a generator
of $\O_\Gamma(1)$ we may identify $V$ with a subspace of
$\O_\Gamma$.  Let $\Gamma'$ be a
subscheme of $\Gamma$, and let $I\subset\O_\Gamma$ be its defining
ideal.  Since $\Gamma$ is in
linearly general position the composite map
$V\rTo \O_\Gamma\rTo\O_\Gamma/I=\O_{\Gamma'}$
is either a monomorphism ($\Gamma'$ is nondegenerate), or an
epimorphism, as one sees directly by comparing the dimension
of the linear space defined by the image of $V$ with the degree of $\Gamma'$.

Now suppose that $\Gamma_1$ is a component of $\Gamma$ of degree $\delta$
which is not Gorenstein, so that $\dim(\socle \O_{\Gamma_1}) >1$.
We will derive a contradiction.
We divide the argument into cases according to the value of $\delta$.

First suppose $\delta\leq r+1$
so that $V\rTo \O_{\Gamma_1}$ is a surjection.  It follows that
the preimage of the socle in
$V$ contains a $\P^{r+2-\delta}$ of hyperplanes,
each of which meets $\Gamma_1$
in a subscheme of degree at least $\delta-1$.
Since $\Gamma$ must have also some
other components, there is a hyperplane in the family meeting
$\Gamma$ in a scheme of degree $\geq \delta-1+(r+2-\delta) = r+1$,
contradicting our assumption of linearly general position.

Suppose now $\delta= r+2$, so that $V\subset \O_{\Gamma_1}$. If
 $\dim(\socle \O_{\Gamma_1}) >1$,
then $V$ meets the socle; thus there is a linear form $x$ on $\P^r$
that meets $\Gamma_1$ in $r+1$ points, again contradicting our hypothesis.

If $\delta= r+3$, then again we have $V\subset \O_{\Gamma_1}$, and we
again get a contradiction as above if $V$ meets the
socle. Thus we may suppose $\O_{\Gamma_1}=V\oplus\socle \O_{\Gamma_1}$,
and that the dimension of the socle is 2.

Set $\m = \m_{\Gamma_1}$. We see that $V\cap \m$
projects isomorphically onto $\m/\m^2$, while
$\socle \O_{\Gamma_1}= \m^2$.
The multiplication on $\O_{\Gamma_1}$
induces a map
$$
V\cap \m \rTo \Hom_k(\m/\m^2, \m^2)
$$
Counting dimensions, and using the algebraic closure of $k$,
we see that for some $x\in V$ the transformation induced by $x$
has rank at most 1. Thus $\O_{\Gamma_1} x$ has dimension at most
2, and  the hyperplane defined by $x=0$ meets
$\Gamma_1$ in at least $r+1$ points, a contradiction.

Finally, suppose that $\delta>r+3$.
By the previous case and \ref{Castel-Gor}, every subscheme of $\Gamma$
of degree $<r+3$ is curvilinear; it follows that every component is
curvilinear, and thus  Gorenstein.
\Box

We prove now a higher dimensional version of
\ref{Castel-Gor}. (See also Cavaliere-Rossi-Valla [1995, Theorem 3.2]
for another proof in a reduced case.)

\theorem{castel-scroll} Let $\Gamma\subset\P^r_k$, $r\ge 3$, be a
finite Gorenstein scheme of degree $\gamma$ which is
in linearly general position, and let $s$ be an integer with $1\le s\le r-2$.
\item{a)} If $\gamma\le r+s+2$, then $\Gamma\subset\P^r_k$ lies on an
$s$-dimensional rational normal scroll  (possibly singular).
\item{b)} If moreover $\gamma\le r+s+2$ and $s\le {{r+1}\over 2}$,
then $\Gamma$ lies on a smooth $s$-dimensional
rational normal scroll in $\P^r_k$.

\proof {\sl a)} The case $s=1$ was proved in part $a)$ 
of \ref{Castel-Gor} so we may assume in the sequel that
$s\ge 2$. It is also enough to prove the assertion when
$\gamma=r+s+2$. Then the Gale transform of $\Gamma$ is
a finite scheme $\Gamma'\subset\P^s_k$ of degree $\gamma$
which is also in linearly general position, by \ref{basic-principle}.
Regarding now $\P^s_k$ as a ``cone'' over $\P^1_k$ with ``vertex''
$\Pi\cong\P^{s-2}_k$ we may resolve $\O_{\Gamma'}(1)$ as
an $\O_{\P^1_k}$-module:
$$0\rTo \oplus_{i=1}^s\O_{\P^1}(-a_i)\rTo
\oplus_{i=1}^{s-1}\O_{\P^1}\oplus\O_{\P^1}(1)\rTo\O_{\Gamma'}(1)\rTo 0,$$
where $\sum_{i=1}^s a_i=\gamma-1=r+s+1$.

Let $\E=\oplus_{i=1}^{s-1}\O_{\P^1}\oplus\O_{\P^1}(1)$, and
let $\F=\oplus_{i=1}^s\O_{\P^1}(-a_i)$. Since $\Gamma'$ is
nondegenerate, there are no sections in $\H^0(\E)=\H^0(\O_{\P^s}(1))$
vanishing identically on $\Gamma'$ and thus $a_i\ge 1$ for all $i$.

To prove the claim we need to check that $a_i\ge 2$ for all $i$, since
then we may use \ref{scroll-gale} to deduce that $\Gamma$ lies on the
birational image of $\P(\F^\ast\otimes\O_{\P^1}(-2))$ in $\P^r$
(since $r+1=\sum_{i=1}^s(a_i-1)$). Twisting by $\O_{\P^1}(1)$ and taking
cohomology in the above short exact sequence we see that
$a_i\ge 2$ for all $i$ iff there are no sections in
$\H^0(\E(1))=\H^0(\oplus_{i=1}^{s-1}\O_{\P^1}(1)\oplus\O_{\P^1}(2))$
vanishing identically on $\Gamma'$. Such sections  correspond to
hyperquadrics in $\P^s$ containing $\Pi$, the ``vertex of the cone'',
and since $\gamma> 2s+1=\h^0(\E(1))$, 
this means that, for a general choice of
$\Pi$, $a_i\ge 2$ for all $i$ iff there are no hyperquadrics
containing both $\Pi$ and $\Gamma$.  The scheme $\Gamma'$ is in
linearly general position so we may conclude by the
following lemma:

\lemma{indep-conditions}
Let $\Gamma\subset\P^s_k$, $s\ge 2$,
be a
finite scheme of degree $\gamma\geq s+3$
which is in linearly general position,
and let $\Pi\subset\P^s_k$ be a general codimension 2 linear
subspace. If $d\leq 3$, then $\Gamma$ imposes $\min(\gamma, ds+1)$ independent
conditions on hypersurfaces of degree $d$ vanishing to order
$d-1$ on $\Pi$.

\conjecture{strong-castel} Castelnuovo's classic result, as generalized to schemes by 
Eisenbud-Harris [1992], says that $\Gamma$ imposes independent 
conditions on forms of degree d. The Lemma above represents a 
strengthening, in that the points impose independent conditions on
a smaller subsystem. We conjecture that the Lemma remains true for
every $d$.

\noindent {\bf Proof of \ref{indep-conditions}} 
(This is \ref{strong-castel} in the case $d\leq 3$.)
Since any scheme of length $\geq s+3$ in linearly 
general position can be extended by the addition of general
points (see Eisenbud-Harris [1992, Theorem 1.3]), we
may assume that $\gamma=ds+1$ and we must show that there are
no $d$-ics  $F$ containing $\Gamma$ and having multiplicity $d-1$ along
$\Pi$. We suppose we have such an $F$ and argue by contradiction.

We first assume $d=2$. Choose a hyperplane $H$ containing a degree $s$
subscheme $\Gamma_1$ of $\Gamma$, and let $\Gamma_2=
\Gamma\setminus\Gamma_1$ be the residual,
a scheme of degree $s+1$.  Specialize $\Pi$ to be a general
$(s-2)$-plane contained in $H$.  
If $F$ contains $H$ and $\Gamma$, then $F=H\cup H'$, where $H'$ is
another hyperplane, containing $\Gamma_2$. This contradicts the
assumption that $\Gamma$ is in linearly general position.  Thus we may
assume that $F$ does not contain $H$, and so $F\cap H=\Pi\cup\Pi'$ for some
$(s-2)$-plane $\Pi'$. As $\Pi$ contains no points of $\Gamma$, $\Pi'$
must contain all of $\Gamma\cap H$, again contradicting the linearly
general position hypothesis.

Now suppose $d=3$. We choose two hyperplanes $H_1,H_2$ as follows:
Order the components of $\Gamma$ by decreasing degree, and suppose the
degrees are $\gamma_1\geq \gamma_2 \geq \dots$. If $\gamma_1\geq 2s$,
then we choose $H_1=H_2$ to be the hyperplane meeting the first
component of $\Gamma$ in $s$ points. Since $\Gamma$ is in linearly
general position and of degree $\geq s+3$ the components of $\Gamma$
are curvilinear, so $H_2$ meets $\Gamma\setminus (H_1\cap \Gamma)$ in
exactly $s$ points as well.

If on the contrary $\gamma_1<2s$, we can divide $\Gamma$ into
two disjoint subschemes each of degree at least $s$ (Reason:
Since no component has degree greater than $2s-1$,
the smallest group of components with total degree $\geq s$ has
total degree $\leq 2s-2$, and the remainder thus has degree $\geq s+3$).
Choose $H_1$ containing a degree $s$ subscheme of one of the 
two subschemes, and $H_2$ containing a degree $s$ subscheme of
the other.

In the first of these two cases, we may choose an $(s-2)$-plane $\Pi$
in $H_1\cap H_2$ that does not meet $\Gamma$; in the second
case the intersection $\Pi=H_1\cap H_2$ automatically misses
$\Gamma$ because $\Gamma$ is in linearly general position.
If $F$ contains both $H_1$ and $H_2$ then after removing the two 
hyperplanes we get a hyperplane containing $s+1$ points of $\Gamma$, a
contradiction as before. If on the other hand $F$ fails to 
contain $H_1$, we restrict to $H_1$. As $2\Pi$ is contained
in $F\cap H_1$, the scheme $H_1\cap\Gamma$ is contained in the
remaining $(s-2)$-plane in $F$ (restricted to $H_1$). 
Once again, this is a contradiction.\Box

\noindent{\sl Proof of \ref{castel-scroll} continued:}
{\sl b)} Again we may assume that $\gamma=r+s+2$, and thus that the
Gale transform $\Gamma'$ is a finite scheme of degree $\gamma$ in $\P^{s}$.
In the notation of {\sl a)}, we need  to check that $a_i\ge 3$ for
all $i$.  Twisting the above short exact sequence by $\O_{\P^1}(2)$,
and taking cohomology this amounts to the fact that no
section in $\H^0(\E(2))$ vanishes identically on $\Gamma'$.
The $3s+1$ sections in $\H^0(\E(2))$ correspond now to
the cubics in $\P^s$ vanishing to second order on $\Pi$, so the claim
follows again from \ref{indep-conditions}.\Box

\remark{} The result in \ref{castel-scroll} is not always sharp, see for
instance \ref{nine3}. 
In case $r\le 2s$, one may slightly improve the statement of
{\sl a)}  by showing that the Gale transform
lies on a rank 4 quadric scroll. \Box

\section{Gale-adjoint} The Gale Transform and Canonical Modules

In this section we consider a fundamental
relation between the presentations of the
canonical modules of a finite Gorenstein subscheme of projective
space and its
Gale transform. 

As in \ref{gale-scheme} we consider a Gorenstein scheme
$\Gamma$ finite of degree $\gamma$ over $k$. Suppose that $\Gamma$ is
embedded in $\P^r_k = \P(V)$ by a linear series $(V,\O_\Gamma(1))$. We
write $W = V^\perp\subset \H^0(K_\Gamma(-1))$ for the linear series
corresponding to the Gale transform.  Supposing that $W$ is
base-point-free, we write $\Gamma^\prime\subset \P^s_k=\P(W)$ for the
image of $\Gamma$ under the corresponding map. Writing $S$ for the
homogeneous coordinate ring of $\P^r$ and $S_\Gamma$ for the
homogeneous coordinate ring of $\Gamma$, we will study the canonical
module $\omega_\Gamma = \Ext^r_S(S_\Gamma, S(-r-1))$.

We will also use the notion of adjoint of a matrix of linear forms.
If $V_1$, $V_2$, and $V_3$ are vector spaces over $k$ and
$\phi\in V_1\tensor V_2\tensor V_3$ is a trilinear form,
then $\phi$ can be regarded as a homomorphism of graded free modules
over the polynomial ring $k[V_1]$:
$$
\phi_{V_1}:\,  V_2^*\tensor k[V_1] \rTo  V_3\tensor k[V_1](1),
$$
and in two other ways corresponding to the permutations of $\{1,2,3\}$.
We call these three linear maps (which may be viewed as matrices of
linear forms, once bases are chosen) {\it adjoints\/} of one another.

\proposition{lin-presentation}  If
$\Gamma^\prime\subset\P(W)$ is the Gale transform
of $\Gamma\subset\P(V)$,
then the linear part of the presentation matrix of
$(\omega_{\Gamma^\prime})_{\geq -1}$, as a $k[W]$-module, is adjoint to
the linear part of the presentation matrix of
$(\omega_\Gamma)_{\geq -1}$ as a $k[V]$-module.

\proof
Consider the multiplication map $V\tensor W\rTo \H^0(K_\Gamma)$,
and let $N$ be its kernel.  From Corollary 1.3 and Theorem 1.4
of Eisenbud-Popescu [1996] we see that $N$ may be regarded
as either the space of linear relations on the degree
$-1$ elements of $\omega_\Gamma$ considered as a module over $k[V]$,
{\it or} as
the space of linear relations on the degree $-1$ elements of
$\omega_{\Gamma^\prime}$ regarded as a module over $k[W]$,
which is exactly the meaning of adjointness.
\Box

To exploit this result we need to know when the linear
part of the presentation matrix of $(\omega_{\Gamma})_{\geq -1}$
actually is the presentation matrix.  This occurs when
$\omega_{\Gamma}$
is generated in degree $\leq -1$ and its relations are generated
in degree $\leq 0$. The first condition is easy to characterize
completely:

\proposition{Koszu-hom} Suppose the field $k$ is
algebraically closed, and let $\Gamma$ be a finite (not necessarily
Gorenstein) scheme in $\P^r=\P^r_k$, not contained in any hyperplane.
Then $\omega_\Gamma$ is generated in degrees $\leq 0$, 
and it fails to be generated in degrees $\leq -1$ iff the homogeneous
ideal $I_\Gamma$ contains (after a possible change of
variables) the ideal of $2\times 2$
minors of a matrix of the form
$$
\pmatrix{
x_0&\dots&x_t&x_{t+1}&\dots&x_r\cr
0&\dots&0&l_{t+1}&\dots&l_r
},
$$
where $0\le t<r$ and the $l_i$ are linearly independent linear forms.

\proof Since $I_\Gamma$ is generated in degree $\geq 2$, its 
$(r-1)^{st}$ syzygies are generated in degree $\geq r$, which yields
the first statement. A standard Koszul homology argument (see
Green [1984] for the source, or for example Cavaliere-Rossi-Valla
[1994], or Eisenbud-Popescu [1997, Proposition 4.2] for this particular
result) shows that the $(r-1)^{st}$ syzygies are generated in degree
$\geq r+1$ unless $I_\Gamma$ contains the ideal of $2\times 2$ minors
of a matrix of the form
$$
\pmatrix{
x_0&\dots&x_t&x_{t+1}&\dots&x_r\cr
l_0&\dots&l_t&l_{t+1}&\dots&l_r
},\leqno{(*)}
$$
where the $l_i$ are linear forms, and the row of $\l_i$ is not a
scalar multiple of the first row.
Because the number of variables is only $r+1$, this $2\times (r+1)$
matrix must have a ``generalized zero'' (Eisenbud [1988]) and
thus may be transformed as in the claim of the proposition.\Box

\corollary{lgp-and-omega} If $\Gamma$ is a finite Gorenstein scheme
of degree $\geq r+2$
in linearly general position in $\P^r$, then $\omega_\Gamma$ is
generated in degrees $\leq -1$.

\proof Suppose not. By \ref{Koszu-hom},
$I_\Gamma$ contains the ideal of minors of a matrix 
with linear entries as in the statement of
\ref{Koszu-hom}.
By our general position hypothesis, the degree of the subscheme
$\Gamma_1$ of $\Gamma$ contained in $V(x_0,\dots,x_t)$ is at most
$r-t$, while the degree of the subscheme $\Gamma_2$ of $\Gamma$
contained in $V(l_{t+1},\dots,l_r)$ is at most $t+1$.  Since
$I_\Gamma$ contains the product of the ideals of these linear spaces,
and $\Gamma$ is Gorenstein, $\deg\Gamma_1+\deg\Gamma_2\geq 
\deg\Gamma$, a contradiction.\Box

In the reduced case, or more generally in the case
when $\Gamma_{\red}$ is nondegenerate, we can give a  geometric
necessary and sufficient condition:

\definition {}
A finite scheme $\Gamma\subset \P^r$ is decomposable
if it can be written as the union of two subschemes contained
in disjoint linear subspaces $L_1$, and  $L_2$, in which case we say that
$\Gamma\subset\P^r$ is the {\it direct sum\/} of its {\it summands\/}
$\Gamma\cap L_1$ and $\Gamma\cap L_2$.

\proposition {decompose} Let $\Gamma\subset\P^r$ be a finite
scheme such that $\Gamma_\red$ is nondegenerate. The module
$\omega_\Gamma$ is generated in degrees $\leq -1$
iff $\Gamma$ is indecomposable.

\proof First suppose that $\omega_\Gamma$
is not generated in degree $\leq -1$. 
By \ref{Koszu-hom}, $I_\Gamma$ contains an ideal
of the form $(x_0,\dots,x_t)\cdot(l_{t+1},\dots,l_r)$, where
the $l_i$ are independent linear forms. If the linear
span of $x_0,\dots,x_t$ in the space of linear forms
meets that of $l_{t+1},\dots,l_r$, then $I_\Gamma$ would
contain the square of a linear form, and thus $\Gamma_\red$
would be degenerate. It follows that
$$(x_0,\dots,x_t, l_{t+1},\dots,l_r)=(x_0,\dots,x_r),$$ so
$$
(x_0,\dots,x_t)\cdot(l_{t+1},\dots,l_r)
=(x_0,\dots,x_t)\cap(l_{t+1},\dots,l_r)
$$
is the ideal of the union of two disjoint linear spaces.

Conversely, if $I_\Gamma$ contains the ideal of the disjoint
union of two linear spaces, then (after a possible change of variables)
it contains an ideal of the form
$$
(x_0,\dots,x_t)\cap(x_{t+1},\dots,x_r)= (x_0,\dots,x_t)\cdot(x_{t+1},\dots,x_r)
$$
which may be written as the ideal of minors of the matrix
$$
\pmatrix{
x_0&\dots&x_t&x_{t+1}&\dots&x_r\cr
0&\dots&0&x_{t+1}&\dots&x_r
},
$$
as required. (The resolution of the product is also easy
to compute directly.)
\Box

The condition that the relations on $\omega_\Gamma$ are generated
in degree 0 is deeper, and is expressed in the third
part of the proposition below (we include the first two parts
because of the nice pattern of results):

\proposition{gen-lin-presentation}
Suppose $\Lambda\subset \P^r_k=\P(V)$ is a finite Gorenstein subscheme
over a field $k$ with algebraic closure ${\overline k}$.
\item {a)} If ${\overline k}\tensor\Lambda$
contains a subscheme of degree $r+1$
in linearly general position, then the $k[V]$-module  $\omega_\Lambda$ is
generated in degree $\leq 0$ and its relations are
generated in degree $\leq 1$.
\item {b)} If ${\overline k}\tensor\Lambda$ contains a subscheme of degree
$r+2$ in linearly general
position, then $\omega_\Lambda$ is
generated in degree $\leq -1$.
\item {c)} If ${\overline k}\tensor\Lambda$
contains a subscheme of degree $r+3$
in linearly general position and ${\overline k}\tensor\Lambda$
does not lie on a
rational normal curve, then the relations on $\omega_\Lambda$
are generated in degree $\leq 0$.
Thus the presentation matrix of ${(\omega_{\Lambda})}_{\geq -1}$
is linear.

\proof The conclusion of each part may be checked after
tensoring with $\overline k$, so we may
assume that $k=\overline k$ from the outset.
The condition of part {\sl a)}
means that the homogeneous ideal $I_\Lambda$ contains
no linear form, and hence $\beta_{i,j}(I_\Lambda) = 0$
for $j<i+2$.  As mentioned above, part {\sl b)} follows from
\ref {Koszu-hom} applied to the subscheme of degree $r+2$,
while  {\sl c)} is
the ``Strong Castelnuovo Lemma'' of Green [1984] and Yanagawa [1994];
see Eisenbud-Popescu [1997]
for a  proof of {\sl c)}
involving syzygy ideals for the Eagon-Northcott
complex. See also Ehbauer [1994]
and Cavaliere-Rossi-Valla [1994] for related results.\Box

Returning to the case of the finite Gorenstein subscheme $\Gamma$ we have:

\corollary{presentation}  Suppose that $\Gamma$ is a  Gorenstein
scheme, finite over a field $k$ with algebraic closure 
$\overline{k}$.  Let $(V,L)$ be a linear series that embeds
$\Gamma$ in $\P^r_k=\P(V)$, and set $(W, K_\Gamma\tensor L^{-1})$ be
the Gale transform.  Let also $s+1=\dim_k(W)$.\hfill\break\indent
If $\overline{k}\otimes 
\Gamma$ contains a subscheme of degree 
$\geq r+3$ in linearly general position and
$\Gamma$ does not lie on a rational normal curve in $\P(V)$, then the
vector space $K := \ker( V\tensor W \rTo \H^0(K_\Gamma))$ has
dimension $rs$, and the corresponding matrix with linear entries 
$K\tensor k[V]\rTo W\tensor k[V](1)$ is a presentation matrix for the
$k[V]$-module ${(\omega_\Gamma)}_{\geq -1}$.

\corollary{can-module}  Suppose that $\Gamma \subset \P^r_k$
is a finite nondegenerate Gorenstein subscheme of
degree $\gamma=r+s+2$, with $r, s\geq 1$, and let $\Gamma'\subset\P^s_k$
be its Gale transform.   The following conditions are
equivalent, and are all satisfied if $\,\overline{k}\otimes\Gamma$
contains a subscheme of degree $r+2$ in linearly general
position:
{\item {a)} $\omega_\Gamma$ is
generated in degree $\leq -1$ as a $k[V]$-module.
\item {b)} $\omega_{\Gamma^\prime}$ is
generated in degree $\leq -1$ as a $k[W]$-module.
\item {c)} The multiplication map
$V\tensor W \rTo \ker(\tau)\subset
\H^0(K_\Gamma)$ is surjective.}\hfill\break
When these conditions are satisfied, both
$(\omega_\Gamma)_{\geq -1}$ and  $(\omega_{\Gamma^\prime})_{\geq -1}$
have precisely $rs$ linearly independent linear relations.

\proof $a)\Rightarrow c)$:
The part of ${(\omega_\Gamma)}_0 = \ker(\tau)$ generated
by $(\omega_\Gamma)_{-1}$ is $V\cdot W$.

\noindent  $c)\Rightarrow a)$: Because $\Gamma$ is nondegenerate,
${\overline k}\tensor \Gamma$ contains a subscheme of length $r+1$
in linearly general position.  By \ref{gen-lin-presentation}, part {\sl a)},
$\omega_\Gamma$ is generated in degree $\leq 0$, so it suffices
to show that $(\omega_\Gamma)_0$ is the image of
$V\tensor (\omega_\Gamma)_{-1} = V\tensor W$.
As $(\omega_\Gamma)_{0} = \ker(\tau)$ we are done.

The symmetry of {\sl c)} completes the proof
of the equivalences.
By \ref{gen-lin-presentation}, the condition of {\sl a)} follows if
$\overline{k}\otimes\Gamma$ contains a subscheme of degree at least $r+2$
in linearly general position.
If the conditions in {\sl a)}, {\sl b)}, {\sl c)} hold, then we can compute
the number of linear relations
in the last statement as $\dim(V)\cdot \dim(W)-
\dim(\ker(\tau)) = (r+1)(s+1)-(r+s+2-1) = rs$.\Box

In the simplest case we can say something about $\omega_\Gamma$
itself:

\corollary{can-lin-present} Suppose that
$\Gamma\subset\P^r_k$ contains a subscheme of degree $r+3$
in linearly general position over $\overline{k}$,
$\Gamma$ imposes independent conditions on quadrics (this occurs
for example when $\gamma \leq 2r+1$ and $\Gamma$ is in
linearly general position)
and that $\Gamma$ does not lie
on a rational normal curve in $\P^r_k$. Then the module $\omega_\Gamma$
has a free presentation by the
$(s+1)\times rs$ matrix of linear forms given in \ref{presentation}.

\proof Since $\Gamma$ imposes independent conditions on quadrics
it is 3-regular, whence
$(\omega_\Gamma)_{-2}=0.$
\Box

In the situation of \ref{can-lin-present}, it is useful to ask about
the adjoint to the presentation matrix of $\omega_\Gamma$, which
is the linear part of the presentation matrix of
$(\omega_{\Gamma^\prime})_{\geq -1}$.  In general we have:

\proposition{flip} Let $\varphi: N\rTo V\tensor W$ be a map of
$k$-vector spaces, with $k$ algebraically closed,
and let $\varphi_V: N\tensor k[V](-1) \rTo W\tensor k[V]$
and $\varphi_W: N\tensor k[W](-1) \rTo V\tensor k[W]$ 
be the corresponding adjoint matrices of linear forms.\hfill\break
\indent  
If the sheafification $L$ of
$\coker(\varphi_V)$ is a line bundle on its support 
$X\subset \P(V)$, then the maximal minors of $\varphi_W$ generate, 
up to radical, the ideal of the image of $X$ under the 
map defined by the linear series $W\subset \H^0(L)$.

\proof Let $\widetilde X$ be the subscheme of
$\P(V)\times \P(W)$ defined by the $(1,1)$-forms in the image of
$\varphi$, and let  $p$ be the projection of
$\P(V)\times \P(W)$ on the first factor. 
On $\P(V)\times\P(W)$ the map $\varphi$ 
corresponds to a morphism $N\tensor\O_{\P(V)\times\P(W)}(-1,-1)
\rTo^\varphi \O_{\P(V)\times\P(W)}$ and  $\varphi_V=p_\ast(\varphi(0,1))$.
Since $R^1p_\ast(\O_{\P(V)\times\P(W)}(-1,0))=0$, it follows
that $p_\ast(\O_{\widetilde X}(0,1))=\coker(\varphi_V)=L$, 
so that the push-forward of a line bundle by $p$ is a line bundle,
which implies that $p_{\mid\widetilde X}$ is an isomorphism
from $\widetilde X$ onto $X$. 

An element of $\P(W)$ is a functional $x:\ W\rTo k$.
It is in the support of $\coker(\varphi_W)$ 
iff $\varphi_W$ drops rank when its entries (elements of $W$)
are replaced by their images under $x$; that is, if
there is a functional
$y:\ V\rTo k$ such that $y\tensor x :\ V\tensor W\rTo k$ annihilates
the image of $N$ via $\varphi$. The projection 
$p_{\mid\widetilde X}$ is an isomorphism
from $\widetilde X$ onto $X$, thus the
projection of $\widetilde X$ to $\P(W)$ 
is the image of the map defined by the sections
$W$ in $L$, and the maximal minors of $\varphi_W$ 
generate, up to radical, the ideal of this image.\Box

In an important special case we can do better:

\corollary{gale-flip} If $\Gamma$ is a finite Gorenstein scheme
and both $\Gamma\subset\P^r$
and its Gale transform $\Gamma'\subset \P^s$ satisfy the hypotheses of
 \ref{can-lin-present}, then $\Gamma^\prime$ has its homogeneous
ideal generated, by the $(r+1)\times(r+1)$ minors of the
adjoint matrix of the presentation matrix of $\omega_\Gamma$.

\proof Set as above $N = \ker(V\tensor W\rTo \H^0(K_\Gamma))$ and
let $\varphi$ denote the inclusion of $N$ in $V\tensor W$;
the matrices $\varphi_V$ and $\varphi_W$ defined in \ref{flip} are, by virtue
of \ref{can-lin-present}, presentations of modules whose sheafifications
are the canonical line bundles on $\Gamma$ and $\Gamma^\prime$,
respectively.  Thus their minors generate the homogeneous ideals of
$\Gamma$ and $\Gamma^\prime$, respectively.
\Box

\example{clebsch} ({\bf The Clebsch transform})
Let $\Gamma\subset\P^2$ be a set of six sufficiently general points. A
familiar transformation in the plane, called the {\it Clebsch
transform\/}, can be constructed from these points: Blow up the six
points, and then blow down the six $(-1)$-curves in the blowup which
are the proper transforms of the 6 conics through five of the six
original points.  The images of the six conics are 6 new points, each
associated to one of the original points (the one through which the
corresponding conic did not pass.)  The new set of six points is the
Clebsch transform of the original set.

Coble [1922] showed that the Clebsch transform of
the six points is the same
as the Gale transform. This follows from
\ref{gale-flip}.
The set $\Gamma$ is cut out by the
maximal minors of a $3\times 4$ matrix $M$
whose cokernel is the canonical module.
By  \ref{gale-flip},
the maximal minors of the $3\times 4$ matrix
$M^\prime$ which is adjoint to $M$ (with respect to the rows)
define the Gale transform $\Gamma^\prime$ of $\Gamma$.
This example also illustrates the case $r=s=2$ of
\ref{determinants-and-Veronese}.

The ``third' adjoint matrix gives the connection of these ideas to the
cubic surface. The determinant of the $3\times 3$ matrix $M^{\prime\prime}$
in four variables which is adjoint to both
$M$ and $M^\prime$ (with respect to their columns)
is the equation of the cubic surface in $\P^3$, image of $\P^2$
via the linear system of cubics through $\Gamma$ (see also
Gimigliano [1989]). The two linear series on the
cubic surface blowing down the two systems of
6 lines described above correspond to the line bundles on the surface
obtained from the cokernel of the $3\times 3$ matrix
(restricted to the surface, where it has constant rank 2)
and from the cokernel of the transpose matrix.
\Box

\section{determinantal} The Gale transform of Determinantal Schemes

Let $\Gamma\subset \P^r$ be a set of points defined by the maximal
minors of a matrix with linear entries vanishing in the generic codimension;
that is,
suppose that there is an $(s+1)\times (r+s)$ matrix $M$ such that the
ideal $I_{s+1}(M)$ generated by the $(s+1)\times(s+1)$ minors of
$M$ defines a set of points ($r=(r+s)-(s+1)+1$). It follows that
the degree of $\Gamma$ is $\gamma={r+s\choose s}$.

In this section
we will see that (in a sufficiently general case) the ${(s-1)}^\st$ Veronese
embedding of this set of points is the Gale transform
of the $(r-1)^\st$ Veronese embedding
of a set of points defined by the adjoint matrix to $M$.
We denote the $d^\th$ Veronese map by $\nu_d:\,\P^r\rTo\P^N$ (where
$N = {r+d\choose r}-1$.)

\theorem{determinants-and-Veronese} Let $V$ and $W$ be $k$-vector
spaces of dimension $r+1$ and $s+1$ respectively.
Let $\phi:\, F\rTo V\tensor W$ be a map of vector spaces with
$\dim_k F = r+s$, and let
$\phi_V:\, F\tensor k[V] \rTo W\tensor k[V](1)$
be the corresponding map of free modules over the
polynomial ring $k[V]$. Let $\phi_W$ be the analogous map over $k[W]$, and
let $\Gamma_V\subset \P^r$ and $\Gamma_W\subset \P^s$
be the schemes defined by the ideals
of minors $I_{s+1}(\phi_V)\subset k[V]$ and
$I_{r+1}(\phi_W)\subset k[W]$, respectively.\hfill\break\indent
If $\Gamma_V$ and $\Gamma_W$ are both zero-dimensional then they are
both Gorenstein, there is a natural isomorphism between them,
and
$$
\nu_{s-1}(\Gamma_V) \quad\hbox{is the Gale transform of}\quad
\nu_{r-1}(\Gamma_W).
$$

The proof will have several steps. We first take care of the
Gorenstein condition:
\proposition {locally-Gorenstein}
With notation as in the theorem, the following are equivalent:
\item{a)} Both $\Gamma_V$ and $\Gamma_W$ are zero-dimensional schemes.
\item{b)} $\codim I_{s+1}(\phi_V) = r$ and  $\codim I_s(\phi_V) = r+1$.
\item{c)} $\Gamma_V$ is zero-dimensional and Gorenstein.

The proof can be analyzed to show that when these conditions are satisfied,
$\Gamma_V$ and $\Gamma_W$ are in fact local complete intersections.

\proof By definition, $\Gamma_V$ is zero-dimensional iff
$\codim I_{s+1}(\phi_V) = r$. Thus to prove the equivalence of
{\sl a)} and {\sl b)} we suppose that
$\Gamma_V$ is zero-dimensional and we must show that
$\Gamma_W$ is zero-dimensional iff $\phi_V$ never drops rank by more than 1.
Let $\phi_F: V^*\tensor k[F^*]\rTo W\tensor k[F^*](1)$
be the third map induced by $\phi$. The
generalized zeros (in the sense of Eisenbud [1988])
of $\phi_F$ in a generalized row indexed by an element
$\alpha\in W^*$ correspond to elements of the kernel of
${\phi_W}_{|\alpha}: F\rTo V$. Thus to say that $\Gamma_W$ is zero-dimensional
is equivalent to saying that only finitely many
generalized rows of $\phi_F$ have generalized zeros. On the other hand
the assumption that $\Gamma_V$ is finite means that only finitely
many generalized columns of $\phi_F$ have generalized zeros, so that
the finiteness of $\Gamma_W$ amounts to saying that no generalized
column can have 2, and thus infinitely many, generalized zeros.  That is,
$\Gamma_W$ is finite iff $\phi_V$ never drops rank
by more than 1.

To prove the equivalence of {\sl b)} and {\sl c)} we may again assume that
$\Gamma_V$ is zero-dimensional.  Since the $(s+1)\times (s+1)$ minors of
$\phi_V$ have generic codimension, we may use the Eagon-Northcott
complex to compute
$$
\omega_{\Gamma_V} =
(\Sym_{r-1}(\coker \phi_V))(s-1).
$$
Now the scheme $\Gamma_V$ is Gorenstein iff $\omega_{\Gamma_V}$
is locally principal. On the other hand,
the ideal $I_{s}(\phi_V)$
has codimension $r+1$ iff it defines the
empty set iff $\coker \phi_V$ is locally principal on $\Gamma_V$
iff
$(\Sym_{r-1}(\coker \phi_V))(s-1)$ is
locally principal on $\Gamma_V$, as required.\Box
\smallskip
\noindent{\sl Proof of \ref{determinants-and-Veronese}}.
We begin with the identification of $\Gamma_V$ and $\Gamma_W$.
Working on $\P := \P(V)\times\P(W)$ the map $\phi$ corresponds to a map
of sheaves
$\phi_{VW}:\, \O_\P(-1,-1)^{r+s}\rTo \O_\P$. We define a subscheme
$\Gamma\subset\P$ by setting $\O_\Gamma := \coker \phi_{VW}$.
Let $p$ be  the projection on the first factor
$p: \P\rTo\P(V)$. We claim that $p$ induces an isomorphism
from $\Gamma$ to $\Gamma_V$. Since the construction is symmetric
in $V$ and $W$, it will follow that $\Gamma$ is naturally isomorphic to
$\Gamma_W$ too, as required.

We have $\phi_V = p_*\phi_{VW}(0,1)$.  Thus
 $p_*\O_\Gamma(0,1) = \coker \phi_V$. By hypothesis $I_s(\phi_V)$
defines the empty set, so $\coker \phi_V$ is a line bundle on
$\Gamma_V$.  In particular, $p(\Gamma)=\Gamma_V$.  By symmetry,
the projection of $\Gamma$ onto the other factor $\P(W)$ is
$\Gamma_W$. Since the fibers of $p$ project isomorphically to
$\P(W)$, and $\Gamma_W$ is zero-dimensional, this shows in particular
that the map $p_{|\Gamma}: \Gamma\rTo \Gamma_V$ is a finite map.
The fact that the push-forward
by $p$ of a line bundle is a line bundle now implies that
$p_{|\Gamma}$ is an isomorphism.

As noted in the proof of \ref{locally-Gorenstein} the presentation
over $k[V]$ of $\omega_{\Gamma_V}$
is given by the Eagon-Northcott complex. It has the form
$$
F\tensor \Sym_{r-2}W\tensor k[V](s-2) \rTo
\Sym_{r-1}W\tensor k[V](s-1)\rTo \omega_{\Gamma_V}\rTo 0,
$$
where the twists by $(s-2)$ and $(s-1)$ indicate as usual
that we regard the first two
terms of this sequence as free modules over $k[V]$ generated in
degrees $(-s+2)$ and $(-s+1)$, respectively. Taking the 
$(s-1)^\st$ Veronese embedding, we see that 
$\omega_{\nu_{s-1}(\Gamma_V)}$ is generated in degree $-1$
with relations generated in degree 0, 
corresponding to the following right-exact sequence:
$$
F\tensor \Sym_{r-2}W\tensor\Sym_{s-2}V
\rTo^\psi
\Sym_{r-1}W\tensor\Sym_{s-1}V
\rTo
(\omega_{\Gamma_V})_0
\rTo
0.
$$
On the other hand, we would obtain the same map $\psi$ if we had
started instead from the scheme $\Gamma_W\subset \P(W)$.  In other
words, the presentation of $\omega_{\nu_{r-1}(\Gamma_W)}$ is adjoint
to the presentation of $\omega_{\nu_{s-1}(\Gamma_V)}$, and by
\ref{lin-presentation} and \ref{gale-flip} 
this means that the two finite schemes are
related by the Gale transform.\Box

\example{ten-det-in-p3}
Let $\Gamma^{\prime\prime}\subset\P^2$ be
a locally Gorenstein scheme of degree ten, not contained in
any plane cubic curve.  It follows from the 
Hilbert-Burch theorem (see for example Eisenbud [1995])
that the ideal of  $\Gamma^{\prime\prime}$ is generated by the
maximal minors of a $4\times 5$ matrix $M''$ with linear entries,
and by \ref{locally-Gorenstein} its 
$2\times 2$ minors generate an irrelevant ideal.
As above, by  \ref{determinants-and-Veronese},
the maximal minors of the $3\times 5$-matrix
$M$ which is adjoint to $M''$ (with respect to the rows)
define a set $\Gamma\subset\P^3$ of ten points,
whose Gale transform $\Gamma^{\prime}\subset\P^5$
coincides with the second Veronese embedding of
$\Gamma^{\prime\prime}$. In this case
the maximal minors of the ``third'' adjoint
$3\times 4$-matrix $M^{\prime}$
in five variables (which is adjoint to both
$M$ and $M^{\prime\prime}$ with respect to their columns)
define a Bordiga surface in $\P^4$, image of $\P^2$
via the linear system of quartics through the
set of ten points $\Gamma^{\prime\prime}$ --- see
Gimigliano [1989] and Room [1938].

A special case  was described
by Coble [1922],  who refers to Conner [1915]
for connections with the geometry of the ``Cayley symmetroid'':
Let $C\subset\P^6$ be a rational
normal sextic curve, and let $X=\Sec(C)\subset\P^6$ be
the secant variety to the curve $C$. $X$ has degree 10,
since this is the number of nodes of a general projection
of $C$ to a plane.  The homogeneous
ideal of $C$ is generated by the $2\times 2$-minors
of either a $3\times 5$ or a $4\times 4$ catalecticant matrix
with linear entries, induced by  splittings of
$\O_{\P^1}(6)$ as a tensor product of two line bundles of
strictly positive degree.
Furthermore, the homogeneous ideal of $X=\Sec(C)$ is generated
by the  $3\times 3$ minors of either of the above two catalecticant
matrices (this is a classical result, see for example 
Gruson-Peskine [1982], or Eisenbud-Koh-Stillman [1988] 
for a modern reference). Let now $\Pi=\P^3\subset\P^6$ be a
general 3-dimensional linear subspace, and
let $\Gamma\subset\P^3$ be a set of ten points
defined by $\Gamma:=\Sec(C)\cap\Pi$.
Then the $2\times 2$ minors of the restriction of the
$3\times 5$ catalecticant matrix generate an irrelevant ideal.
So as above, by \ref{determinants-and-Veronese},
$\Gamma'$ the Gale transform  of $\Gamma$ lies on
a quadratic Veronese surface in $\P^5$.

In this case, the maximal minors of the ``third'' adjoint matrix
cut out a special Bordiga surface in $\P^4$, image of $\P^2$
via the linear system of quartics through the
ten nodes of the rational plane sextic curve obtained by
projecting the rational normal curve above from the $\P^3$.
(See  also Room [1938], 14.21, p.~391 and ff,
Hulek-Okonek-van de Ven [1985], and Rathmann [1989].)

\problem{} Can the previous method be used to tell exactly
when a set of 10 points in $\P^3$ is  determinantal?
\medskip

For geometry related to the Gale transform of
$\gamma\ge 11$ points on the Veronese
surface in $\P^5$  see also Davide [1997].

\medskip

\section{Gorenstein} Gorenstein and Self-associated Schemes

Castelnuovo's and Coble's original interest in associated sets of points
centered on those sets that are ``self-associated'', that is
equal to their own Gale transform up to projective equivalence.

Dolgachev and Ortland [1988] posed the problem of giving a ``clear-cut
geometrical statement'' equivalent to self-association (Remark 3,
p. 47). The following is our solution to this problem:

\theorem{self-assoc-crit} Let $\Gamma\subset\P^r=\P(V)$ be a finite
Gorenstein scheme of degree $2r+2$ over an algebraically closed
field $k$. The following are equivalent:
\item{a)} $\Gamma$ is self-associated.
\item{b)} Each of the (finitely many) 
subschemes of degree $2r+1$ of $\Gamma$
imposes the same number of conditions on quadrics as $\Gamma$ does.
\item{c)} If we choose a generator of $\O_\Gamma(1)$, and thus
identify $V$ with a subspace of $\O_\Gamma$, 
there is a linear form $\phi:\O_\Gamma\rTo k$ which vanishes
on $V^2$ and which generates $\Hom_k(\O_\Gamma,k)$ as an
$\O_\Gamma$-module.

We include $c)$ because it represents the
most efficient way we know to check the property of self-association
computationally. Namely, representing the multiplication
table of the ring $\O_\Gamma$ as a matrix with linear entries
over $\Sym(\O_\Gamma)$, we may identify $V^2$ with a vector
space of linear forms in $\Sym(\O_\Gamma)$. Then part
$c)$ in \ref{self-assoc-crit} can be reformulated as:
$\Gamma$ is self-associated iff the matrix of the
multiplication table reduced modulo the linear forms in $V^2$
has maximal rank. This test can be implemented
in Macaulay/Macaulay2.

\proof  The subscheme $\Gamma$ is self-associated iff there is
an isomorphism of $\O_\Gamma$-modules
$\varphi:\,\O_\Gamma(1)\rTo\O_\Gamma(1)^*$
such that the composite of natural maps
$$
V\rTo \H^0(\O_\Gamma(1))\rTo^\varphi \H^0(\O_\Gamma(1)^*)\rTo V^*
$$
is zero. 

Giving a morphism of $\O_\Gamma$-modules
$\varphi:\,\O_\Gamma(1)\rTo\O_\Gamma(1)^*$
is the same as giving a map of vector spaces
$\overline\varphi:\,\O_\Gamma(2)=\O_\Gamma(1)\tensor\O_\Gamma(1)\rTo k$.
Since $\O_\Gamma$ is Gorenstein, the modules
$\O_\Gamma(1)$ and $\O_\Gamma(1)^*$ are isomorphic, and
$\varphi$ is an isomorphism iff 
$\overline\varphi$ generates $\Hom_k(\O_\Gamma(2),k)$ as an $\O_\Gamma$-module.
We may write $\O_\Gamma=\prod\O_{\Gamma_i}$, where the $\Gamma_i$
are the connected components of $\Gamma$, and the ideals
$J_i=\socle(\O_{\Gamma_i})$ are all 1-dimensional. With this
notation, $\varphi$ is an isomorphism iff 
$\overline\varphi$ does not annihilate any of the one-dimensional
submodules $J_i\O_\Gamma(2)$.

On the other hand, $\varphi$ makes the composite map displayed above zero
iff $\overline\varphi$ annihilates $V^2$, the image of $V\tensor V$ in
$\O_\Gamma(2)$. Thus $\Gamma$ is self-associated iff there is a map
$\overline\varphi$ that annihilates $V^2$ but not any of the $J_i\O_\Gamma(2)$
 iff $V^2\cap J_i\O_\Gamma(2)=0$ iff each codegree 1 subscheme of $\Gamma$
imposes the same number of conditions on quadrics as $\Gamma$, 
thus proving that $a)$ and $b)$ are equivalent.

Once we choose an  identification
of $\O_\Gamma$ and $\O_\Gamma(1)$, part $c)$ is a reformulation of 
this condition.\Box

A classical theorem of Pascal says that given a conic in the plane and
two triangles circumscribing it, (algebraically this means that the
vertices of each triangle are apolar to the conic), then the six
vertices of the two triangles all lie on another conic.  In other
words the six points form a set of self-associated points in the
plane, as one sees from \ref{self-assoc-crit}.

Coble [1929] generalized this statement to say that
for sufficiently general sets of
$2r+2$ points in $\P^r$, self-association is the same as
failing to impose independent conditions on quadrics.  
We next characterize arithmetically Gorenstein schemes
in terms of the Gale transform, and we will see in 
a somewhat more precise way
that a self-associated
scheme in $\P^r$ is the same as an arithmetically Gorenstein scheme of 
degree $2r+2$ except in degenerate circumstances:

\theorem{vero-gale} If $\Gamma\subset\P^r$ is a finite nondegenerate 
Gorenstein scheme, then $\Gamma$ is arithmetically Gorenstein iff
$\omega_\Gamma$ is generated in degrees $\leq -1$ and the Gale
transform of $\Gamma$ is the $d^\th$ Veronese embedding of $\Gamma$
for some $d\geq 0$.  In particular if $\deg(\Gamma)=2r+2$, then
$\Gamma$ is arithmetically Gorenstein iff $\omega_\Gamma$ is generated
in degrees $\leq -1$ and $\Gamma$ is self-associated.

The case $d=0$ occurs only for $r+2$ points in $\P^r$; such a 
scheme is arithmetically Gorenstein iff it is in linearly 
general position.

\proof If $\Gamma$ is nondegenerate and arithmetically Gorenstein,
then the symmetry of the free resolution of the homogeneous
coordinate ring $S_\Gamma$ shows that $\omega_\Gamma=S_\Gamma(d+1)$
for some $d\geq 0$, and $\omega_\Gamma$ is generated in degree
$-d-1\leq -1$. By \ref{hypersurface-series} the Gale transform
is given by the image of $(\omega_\Gamma)_{-1}={(S_\Gamma)}_d$
in $\H^0(K_\Gamma(-1))$, so the Gale transform is the 
$d^\th$ Veronese embedding.

Conversely, suppose $\omega_\Gamma$ is generated in degrees $\leq -1$
and the Gale transform coincides with the $d^\th$ Veronese 
embedding for some $d\geq 0$.
Since $S_\Gamma(d+1)$ is also generated in degrees $\leq -1$,
and both $\omega_\Gamma$ and $S_\Gamma(d+1)$ are Cohen-Macaulay
modules, they are isomorphic iff
${(\omega_\Gamma)}_{\geq -1}\iso {(S_\Gamma(d+1))}_{\geq -1}$, and
this occurs precisely when there is an isomorphism of sheaves of
$\O_\Gamma$-modules $K_\Gamma(-1)\iso \O_\Gamma(d)$, which maps
$V^\perp$ to ${(S_\Gamma)}_d$. This last is the condition that the
Gale transform of $\Gamma$ is the $d^\th$ Veronese embedding
of $\Gamma$.\Box

Perhaps the best characterization of this kind is

\theorem{fail} If $\Gamma\subset\P_k^r$ is a nondegenerate finite scheme
of degree $2r+2$ over an algebraically closed field $k$,
then $S_\Gamma$ is Gorenstein if and only if $\Gamma$ is
self-associated and fails by 1 to impose independent conditions
on quadrics.

By \ref{self-assoc-crit}, we could restate the condition of the Theorem by
saying that $\Gamma$ fails to impose independent 
conditions on quadrics but that every maximal proper subscheme of $\Gamma$ 
(equivalently every proper subscheme of $\Gamma$)
imposes independent conditions on quadrics.
This is also a consequence of Kreuzer [1992], Theorem 1.1,  
which generalizes the main result of 
Davis-Geramita-Orecchia [1985] to the non-reduced case.

A result of Dolgachev-Ortland [1988] (Lemma 3, p. 45) and
Shokurov [1971] shows that every proper subscheme 
does impose independent conditions if $\Gamma$ is reduced and, 
for every $s<r$, no subset of $2s+2$ points of $\Gamma$ is 
contained in a $\P^s$, which is the same as
saying that $\Gamma$ is {\it stable\/} in this case.
(See \ref{stability} for the general stability test.) 
Dolgachev and Ortland [1988]  use
this to prove that a reduced set of stable points is self-associated
if and only if it fails to impose independent conditions on quadrics,
generalizing a result of Coble [1929].
\medskip

\noindent{\sl Proof of \ref{fail}}.
If $S_\Gamma$ is Gorenstein then $\Gamma$ is self-associated by
\ref{vero-gale}, and fails by just 1 to impose independent conditions
on quadrics since
$(S_\Gamma)_0=(\omega_\Gamma)_{-2}=((S_\Gamma)_2)^\perp$
is 1-dimensional.

Conversely, suppose that $\Gamma$ is self-associated. 
By definition there is
an isomorphism of $\O_\Gamma$-modules 
$\O_\Gamma(1)\rTo K_\Gamma(-1)$ carrying
$V={(S_\Gamma)}_1$ to ${(\omega_\Gamma)}_{-1}$.  This defines 
a map of modules ${(S_\Gamma(2))}_{\geq 1}\rTo \omega_\Gamma$. Since
$\Ext^1_{S_\Gamma}(k,\omega_\Gamma)=k$, concentrated in degree 0, this
map lifts to a map 
$\alpha:\,(S_\Gamma(2))\rTo \omega_\Gamma$, 
necessarily a monomorphism.  As $\Gamma$ is
nondegenerate, $\omega_\Gamma$ is generated in degree $\leq 0$. We
have ${(\omega_\Gamma)}_0=({(S_\Gamma)}_0)^\perp$, so it has dimension
just one less than $\deg \Gamma$.  

If further $\Gamma$ fails by 1 to
impose independent conditions on quadrics, then $(S_\Gamma)_2$ has the
same dimension, and we see that $\alpha$ is an isomorphism in all
degrees $\geq -2$. Further, if every subscheme of  $\Gamma$
imposes independent conditions on quadrics then $\Gamma$ 
imposes independent conditions on cubics (Proof: Find a quadric
vanishing on a codegree 2 subscheme. It does not 
generate a minimal submodule of $\O_\Gamma(2)$, so
we can multiply by a linear form to get a cubic vanishing precisely
on a codegree 1 subscheme.) Thus $(\omega_\Gamma)_{-d}=0$
for $d\geq 3$, and $\alpha$ is an isomorphism. Thus
$\Gamma$ is arithmetically Gorenstein.\Box

As mentioned above, if we restrict to the case of stable sets of 
points, there is a particularly simple characterization of 
self-association, due to Coble [1929].

\corollary{quadric-fail} (Coble, Dolgachev-Ortland)
A stable set of $2r+2$ distinct $k$-rational points in $\P_k^r$ is
self-associated iff it fails to impose independent conditions on quadrics.

\proof Any self-associated scheme of degree $2r+2$
fails to impose independent
conditions on quadrics, for example by 
\ref{self-assoc-crit}.
Conversely, assume that $\Gamma$ is stable and fails to impose independent
conditions on quadrics. Any
subset of $2r+1$ points of $\Gamma$ 
imposes independent conditions on quadrics
(see Dolgachev-Ortland [1986], Lemma 3, p.45 and Shokurov [1971]),
and so the result follows from
\ref{self-assoc-crit}.\Box

As a corollary of \ref{fail} we can exhibit an interesting class of
examples. To
simplify the notation we systematically identify effective
divisors on a smooth
curve with the schemes they represent.

\corollary{K+x+y} Let $C$ be a reduced irreducible canonically embedded
 curve in $\P^n$,
and let $\Gamma\subset C$ be a Cartier divisor in the class
$K_C+D$, where $D$ is an effective divisor of degree 2,
so that the degree of $\Gamma$ is $2n+2$.
The scheme $\Gamma$ is arithmetically Gorenstein iff $\Gamma$
does not contain $D$.

\proof We check the conditions of \ref{fail}, noting that by Riemann-Roch
an effective  Cartier divisor $E$ fails to impose independent conditions on
quadrics iff 
$$
\h^1(2K_C-E)=\h^0(E-K_C)\neq 0.
$$ 
In particular, any effective
divisor in the
class $K_C+D$ fails to impose independent conditions on quadrics
because $D$ 
imposes just one condition. Further, if $\Gamma\supset D$
and $p\in D$, then $\h^0(\Gamma-p-K_C)=\h^0(D-p)\neq 0$, and
we see that $\Gamma$ is not arithmetically Gorenstein.

Now suppose that $\Gamma$ does not contain $D$, that is, 
$\Gamma-D$ is ineffective. As $\h^0(K_C) = \h^0(K_C+D)-1$, we
see $\Gamma$ cannot even contain a point of $D$. Thus, 
for any $p\in\Gamma$, $\h^0(\Gamma-p-K_C)=\h^0(D-p)=0$,
so $\Gamma-p$ imposes independent conditions on quadrics,
and $\Gamma$ is arithmetically Gorenstein by \ref{fail}.\Box

We conclude with a remark that will be used for the
classifications in \ref{classification}:
Self-associated schemes can be
direct sums (in the sense of
\ref{Gale-adjoint} above):

\proposition {degen-self-assoc} If $\Gamma\subset\P^r$ is a
decomposable finite scheme, then $\Gamma$ is self-associated
iff each of its summands is self-associated.

\proof If $\Gamma$ decomposes with summands
$\Gamma_1\subset L_1$ and $\Gamma_2\subset L_2$,
then the natural map $V\rTo\H^0(\O_\Gamma(1))$ splits as the direct sum
of maps $V_i\rTo\H^0(\O_{\Gamma_i}(1))$, where $V_i=\H^0(\O_{L_i}(1))$.
Thus the Gale transform $V^\perp\subset \H^0(K_\Gamma(-1))$ splits
too, and the proposition follows.\Box

\corollary{agorenstein}  A finite locally Gorenstein
scheme $\Gamma\subset\P^r$
of degree $2r+2$  which is
a direct sum of arithmetically Gorenstein
subschemes is self-associated. The number of summands is exactly
the amount by which $\Gamma$ fails
to impose independent conditions on quadrics.\Box

\section{linalg} Linear Algebra and Self-association

Let $\Gamma\subset \P^r_k=\P(V)$ be  a nondegenerate set of $2r+2$ distinct
points.
Choose  $\Gamma_1\subset \Gamma$
a subset of $r+1$ points that spans $\P^r$, and
let $\Gamma_2$ be the complementary set.  If $\Gamma$ is self-associated,
then by \ref{basic-principle} the set $\Gamma_2$ must also span $\P^r$.

Babbage [1948] pointed out that $\Gamma$ is self-associated iff in
addition there is a nonsingular quadric $Q\subset \P^r$ such that each
of $\Gamma_1$ and $\Gamma_2$ are {\it apolar} (self-conjugated
simplexes) with respect to $Q$. In modern language (and replacing
quadratic forms with symmetric bilinear forms to avoid problems in
characteristic 2):

\theorem{bath-orth-bases}
Let $\Gamma\subset \P^r_k$ be  a nondegenerate set of $2r+2$ distinct points.
The set $\Gamma$ is self-associated if and only if it can
be decomposed into a disjoint union $\Gamma=\Gamma_1\cup\Gamma_2$,
where $\Gamma_1$ and $\Gamma_2$ correspond to orthogonal bases for
the same nonsingular symmetric bilinear form on $V$.

\proof The result follows immediately from Castelnuovo's
definition (in \ref{history}) which amounts to saying that set of
points is self-associated iff there is a nonsingular bilinear form
$B$ on
$V$ and a decomposition of $\Gamma$ into two disjoint bases $\{e_i\}$ and
$\{f_i\}$ for $V$ (the vertices of the two simplices) such that
the orthogonal complement of
any $e_i$ is the span of $\{e_j\}_{j\neq i}$,
 and similarly for the $f_i$.  Since
the form $B$ has orthogonal bases, it is symmetric. \Box

The bilinear form really does depend on the choice of splitting
$\Gamma=\Gamma_1\cup\Gamma_2$, since otherwise each vector in $\Gamma$
would be orthogonal to all the other vectors in $\Gamma$, and these
other vectors span $V$.
Babbage also asserts that $Q$ is unique, which, as we shall see, is
false in general.  

To generalize Babbage's result to schemes, we first extend the
notion of an orthogonal basis.

\definition{orthogonal-basis} As above, let
$V$ be a $k$-vector space of
dimension  $\dim_k(V) =r+1$.
\item{a)} A scheme $\Lambda\subset \P(V)$ is a {\it basis\/} for
$V$ if it is nondegenerate and of degree $r+1$ over $k$;
that is, if the natural map $V\rTo \H^0(\O_\Lambda(1))=\O_\Lambda(1)$
is an isomorphism.
\item{b)} If $B:\; V \rTo V^*$ is a $k$-linear map,
then we say $\Lambda$ is an {\it orthogonal basis\/}
with respect to the bilinear form
(or quadratic form) corresponding to $B$
if $\Lambda$
is a basis and the composition
$$
\O_\Lambda(1)\iso V\rTo^B V^*\iso K_\Lambda(-1)
$$
is an isomorphism of sheaves of $\O_\Lambda$-modules.

This generalizes the
classical notion: If $\Lambda$ corresponds
to an ordinary basis $\{p_i\}_{\{i=1,\ldots,r+1\}}$ of $V$, then
$\O_\Lambda \iso k\times k\times\dots\times k$ as a ring,
and $\O_\Lambda$ has
idempotents $e_i$ corresponding to the basis elements $p_i$.
Any sheaf $\F$ over $\O_\Lambda$ decomposes as $\oplus e_i\F$,
and any morphism of sheaves preserves this decomposition.  Thus
$B$ satisfies condition $b)$ above
iff the matrix of $B$ with respect to the basis $\{p_i\}$ is
diagonal, and $B(p_i)(p_j)=0$ for all $i\neq j$.

In the classical case the existence of an orthogonal basis
implies that the bilinear form corresponding to $B$ is
symmetric. This remains also true in our generality:
For if $B$ is a sheaf homomorphism, then
for any generator $f\in \O_\Lambda(1)$ and section
$g\in K_\Lambda(-1)$ the quotient is a section $g/f\in \O_\Lambda$,
and we have
$$
B(f)(g) = B(f)((g/f)f) = B((g/f)f)(f) = B(g)(f).
$$
Since any section $f'\in \O_\Lambda(1)$ may be written as
$f'=rf$ for some  $r\in \O_\Lambda$, we get
$$
B(f')(g)= rB(f)(g)=rB(g)(f)=B(g)(rf)=B(g)(f'),
$$
as required.

We will generalize  \ref{bath-orth-bases} as follows:

\theorem{orth-bases} Suppose that a finite Gorenstein scheme
$\Gamma\subset\P^r_k=\P(V)$
of degree $2r+2$ decomposes as the disjoint union of two
subschemes $\Gamma_1$ and $\Gamma_2$ that are bases.
Then the scheme $\Gamma$ is self-associated iff
$\Gamma_1$ and $\Gamma_2$
are both orthogonal bases for the same nonsingular bilinear form
on $V$.  The bilinear form is unique iff $\Gamma$ is arithmetically Gorenstein.

\proof Suppose first that $\Gamma$ is
self-associated, so there is a sheaf
isomorphism $\lambda: \O_\Gamma(1)\rTo K_\Gamma(-1)$ carrying
$V$ to $V^\perp$.
As $\lambda$ is a sheaf homomorphism, it decomposes as
a direct sum of isomorphisms
$$
\lambda_i:\; \O_{\Gamma_i}(1)\rTo K_{\Gamma_i}(-1), \quad i=1,2.
$$
Write $a_i:\; V\rTo \O_{\Gamma_i}(1)$ for the
inclusion corresponding to the embedding of $\Gamma_i$ in
$\P(V)$.
The bilinear forms $B_i=a_i^*\circ\lambda_i\circ a_i:\; V\rTo V^*$, $i=1,2$,
are nonsingular since each of the three maps in the composition
is an isomorphism. Thus each $\Gamma_i$ is an orthogonal
basis for $B_i$.  As $\lambda = \lambda_1\oplus\lambda_2$
maps $V$ to $V^\perp$, the bilinear form
$$
\pmatrix{a_1^* & a_2^*}\circ
\pmatrix{\lambda_1&0\cr 0&\lambda_2\cr}\circ
\pmatrix{a_1 \cr a_2}
$$
is zero. Thus $B_1=-B_2$, and each $\Gamma_i$ is an orthogonal basis
for $B_1$. 

The bilinear symmetric forms $B$ making $V$ isotropic 
correspond to
elements of the dual of the cokernel of 
$\Sym_2(V)\rTo\O_\Gamma(2)$.  If $\Gamma$
is self-associated then, by \ref{fail}, $\Gamma$ is
arithmetically Gorenstein iff it fails by exactly one to impose
independent conditions on quadrics. Thus the bilinear form is uniquely
determined up to a scalar factor exactly in this case.

Conversely, if both $\Gamma_i$ are orthogonal bases for a nonsingular
form $B$, then $B$ induces an isomorphism of $\O_\Gamma$-modules
$$
\pmatrix{B& 0\cr 0& -B\cr}
:\; \O_\Gamma(1)=
\oplus_{i=1,2}\O_{\Gamma_i}(1)
\rTo
\oplus_{i=1,2}K_{\Gamma_i}(-1)=
K_\Gamma(-1),
$$
whose associated bilinear form $\lambda$ satisfies the conditions
of \ref{self-assoc-crit}.\Box

As a consequence of \ref{orth-bases} we can give a new proof of the following
result:

\corollary{flag-var} (Coble, Dolgachev-Ortland)
The variety of ordered arithmetically Gorenstein sets of $2r+2$ distinct
$k$-rational points in $\P^r_k$ whose first $r+1$ elements span,
up to projective equivalence,
is isomorphic to an open subset in the variety of complete flags in
$\P^r_k$. In particular it is irreducible and rational of dimension
${r+1\choose 2}$, and thus the variety of unordered self-associated sets of
distinct points is also irreducible and unirational of the same dimension.

\remark{} In the case of 6 points in the plane, the unordered self-associated
sets of 6 points form a rational variety, isomorphic to the moduli space of
genus 2 curves (Igusa [1964]): such a set of points lies on a conic, whose
double
cover branched over the 6 points is a curve of genus 2.  Is the
variety of unordered self-associated sets always rational?

\medskip

\noindent{\sl Proof of \ref{flag-var}.} If $\Gamma=\Gamma_1\cup \Gamma_2$ is an
arithmetically  Gorenstein set in $\P^r_k$ decomposed into its subsets of the
first $r+1$ and last $r+1$ points,
and $\Gamma_1$ spans $\P^r_k$, then, by \ref{basic-principle},
the set $\Gamma_2$ also spans. By \ref{orth-bases} there is a
unique symmetric bilinear form $B$ for which both $\Gamma_1$ and $\Gamma_2$ are
orthogonal bases. There is a unique
projective equivalence taking $\Gamma_1$ to the
standard simplex of $\P^r_k$ and taking $B$ to the form whose matrix with
respect to this
basis is the identity matrix. To the set $\Gamma$ we may associate the flag
consisting of
the spaces spanned by the first $i$ elements of $\Gamma_2$, for all
$i=1,\dots,r+1$.

Conversely, let $\Gamma_1$ be a basis in $V$,  
and let $B$ be the bilinear form
whose matrix with respect to this basis is the identity. 
Suppose that $V_1\subset\dots\subset V_{r+1}$ is a flag with
$\dim_k(V_i)=i$ transverse to the flag of coordinate subspaces
defined by $\Gamma_1$.
Assume that the
restriction
of $B$ to each subspace $V_i$ is nonsingular. By the Gram-Schmidt process, we
may choose an orthogonal  basis
$v_1,\dots, v_{r+1}$ for $B$ such that $V_i= <v_1,\dots,v_i>$,
$i\in\{1,\dots,r+1\}$.
For an open set of
flags, the points in $\P^r$ corresponding to 
the vectors $v_i$ will be distinct from
the points of the standard simplex. Let $\Gamma_2$ be the set of 
these points, and
set $\Gamma=\Gamma_1\cup \Gamma_2$. By
\ref{orth-bases}, the set $\Gamma$ is arithmetically Gorenstein.\Box

By using \ref{orth-bases} we may partially decide when is it
possible to extend a finite set (or more generally a locally Gorenstein finite
scheme)
to an arithmetically Gorenstein one.

\theorem{complete}
A general set $\Gamma\subset\P^r_k$ of $\gamma=r+1+d\leq 2r+2$ points
can be extended to an arithmetically Gorenstein set $\Gamma\cup
\Gamma'\subset\P^r_k$ of $2r+2$ points iff ${d\choose 2}\leq r$.
Moreover in case ${d\choose 2}=r$, there is a unique linear subspace
$L\subset\P^r$ of dimension $(r-d)$ such that if
$\Gamma\cup\Gamma'$ is arithmetically Gorenstein then
$\Gamma'$ spans $L$.

\proof Changing coordinates if necessary, we may assume without
loss of generality that $\Gamma$ contains as a subset the standard basis
$\Gamma_1=\{e_0,e_1,\dots,e_r\}$ of the ambient vector space $V$.
We write $\Gamma = \Gamma_1\cup\Sigma$, where $\Sigma$ consists
of the remaining $d$ points.

To find an arithmetically Gorenstein scheme
$\Gamma_1\cup\Sigma\cup\Sigma'$ we search for bilinear forms
$B$ such that both $\Gamma_1$ and $\Sigma\cup\Sigma'$ are
orthogonal bases. First, $\Gamma_1$ is an orthogonal basis iff
the matrix of $B$ is diagonal. The mutual orthogonality of 
the elements of $\Sigma$ imposes ${d\choose 2}$ homogeneous
linear equations
on the $r+1$ diagonal elements of $B$.
In particular, the
system has a non-trivial solution whenever ${d\choose 2}\le r$.
Since $\Gamma$ is general and self-associated sets do exist,
there is a solution $B$ which is nonsingular. Choosing $\Gamma_2$ to
be any orthogonal basis containing $\Sigma$, and using \ref{orth-bases}, we
prove
the first statement of the theorem.

To complete the proof we show that for a general set of
points $\Gamma$ the ${d\choose 2}$
linear equations on the coefficients of $B$ above are of maximal rank.
It follows that in case $r+1={d\choose 2}$ the form $B$ is 
unique, and $\Sigma'$ spans the orthogonal complement of $\Sigma$
with respect to $B$. On the other hand, if $r+1<{d\choose 2}$ 
then the equations
for $B$ have only the trivial solution.

Consider the incidence variety $\I$ whose points are pairs
consisting of a nonsingular
bilinear form $B$ on $V$ such that $\Gamma_1$ is an orthogonal basis,
and a $d$-tuple of distinct points in $\P^r$ which
are non-isotropic and mutually orthogonal with respect to the bilinear form
$B$.
To show that the set of linear equations above is linearly independent, we must
show
that the fiber of $\I$ over a general set $\Sigma$ has the expected dimension,
which is  $r-{d\choose 2}.$ We know from the argument above that the general
fiber  has dimension at least the expected dimension. Further, the
set of $d$-tuples of points $\Sigma$ has dimension $rd$. Thus it suffices to
show that $\dim\I = r(d+1)-{d\choose 2}$.

Projecting a pair $(B,\Sigma)\in \I$ onto the first factor we obtain a
surjection $\I \rTo^{} (k^*)^r$. The fiber over a point $B\in {(k^*)}^r$
may be identified with the set of flags $V_1\subset\dots\subset V_d\subset V$,
where $\dim_k(V_i) = i$
and $B$ restricted to each $V_i$ is nonsingular:
given any $\Sigma = \{v_1,\dots,v_d\}$ we let $V_i=\langle
v_1,\dots,v_i\rangle$, and
conversely given the flag we use the Gram-Schmidt process to produce an
orthogonal
basis.  Thus the fiber is an open set in a flag variety of dimension
$rd-{d\choose 2}$,
and so $\dim\I = r(d+1)-{d\choose 2}$ as required. (One can show
further that the
incidence variety is irreducible and nonsingular, but we do not need this.)
\Box

\example{} Five general points in $\P^2$ lie on a unique conic.
Any sixth point on the conic gives an arithmetically Gorenstein set.

\example{} Let $\Gamma$ be a set of seven points in linearly
general position in $\P^3$. They lie on just three independent
quadrics. If these form a complete intersection, then there exists
a unique extension of the seven points to a 
self-associated scheme of degree 8.
If not, the three quadrics cut out a twisted cubic curve
and $\Gamma$ lies on it.  In this case there are many
possible extensions: we can add any further point on the rational normal
curve, and these are the only possibilities.

\example{} Consider now a set $\Gamma\subset\P^6$ of
11 general points. By \ref{complete}, the set $\Gamma$
may be completed to a self-associated set of 14 points
in $\P^6$. For all possible completions the linear
span of the extra three points is a distinguished plane $\Pi=\P^2\subset\P^6$.
The Koszul complex built on the equations defining this
2-dimensional linear subspace
is the complex $E_\bullet(\mu)^\ast(-8)$ embedded at the back end
of the minimal free resolution of the homogeneous ideal $I_\Gamma$ (see the
end of the introduction of Eisenbud-Popescu [1996] for notation and
details). A similar remark holds for a general set $\Gamma$ of
${s+2\choose 2}+1$ points in $\P^{{s+2\choose 2}-s-1}$, yielding a
distinguished linear subspace of dimension ${s+2\choose 2}-2s-2$ whose
equations contribute to the resolution of $\Gamma$ in an interesting way.

\medskip
Since we have been dealing with the condition of forming two orthogonal
bases, we comment on the condition that a set of $2r+2$ points in $\P^r_k$
correspond
to the union of 2 bases. First recall the criteria of stability and
semistability:

\proposition{stability} (Dolgachev-Ortland). Let
$\Gamma \subset \P^r_k$ be a set of $\gamma$ points. Then
$\Gamma$ is semistable if and only if for all $m$ with
$1\le m\le \gamma-1$ the projective linear
span of any subset of $m$ points of $\Gamma$  has dimension
at least ${m(r+1)/\gamma}-1$. Similarly,
$\Gamma$ is a stable set of points if and only if all
previous inequalities  are strict. \Box

The special case when
$\gamma=2r+2$ has a nice linear algebra interpretation:

\lemma{two-bases} (Edmonds [1965]; see also Eisenbud-Koh [1987])
A set $\Gamma$ of $2r+2$ points in $\P^r$ is semistable
iff the points of $\Gamma$ form two bases for the underlying vector space of
the
ambient projective space.
\Box

By the remarks at the beginning of this section, any self-associated set
is semistable.

\section{classification} Classification of Self-associated Schemes in Small
Projective Spaces

In this section we give a complete classification of 
self-associated schemes in $\P^r$ for $r\leq 3$,
and we review classification results of
Coble, Bath, and Babbage for $\P^4$ and $\P^5$. 
We begin with some examples valid in all dimensions,
coming from 
\ref{Goppa} and \ref{K+x+y}.

\proposition{self-examp} The following are families of arithmetically
Gorenstein
nondegenerate schemes of degree $2r+2$  in $\P^r$:
\item{a)} A Cartier divisor in the class $2H-K_C$ on a rational normal
curve $C\subset\P^r$ (defined by the minors of a $2\times r$ matrix  with linear entries).
\item{b)} A quadric section of a nondegenerate reduced irreducible
curve of degree $r+1$, and arithmetic genus 1 in $\P^{r}$, $r\ge 2$.
\item{c)} A hyperplane section of a canonical curve of genus $r+2$, $r\ge 1$.
\item{d)} A Cartier divisor $\Gamma$ in the class $K_C+D$ on a curve
$C$ of genus $g=r+1$, $r\ge 2$, where $D$ is effective of degree 2,
and $\Gamma$ doesn't contain $D$.\Box

The families listed in \ref{self-examp} account in fact
for the general self-associated sets of points in small
projective spaces, as we will see bellow. An easy count of
parameters shows that in \ref{self-examp}, the families
described in $a)$ and $b)$  have dimensions
$2r-1$, $2r+2$ when $r\ge 4$, and 6 when $r=3$, respectively.
For the last claim we use:

\proposition{} Let $\Gamma\subset\P^r$ be a quadric section of an elliptic
normal curve in $\P^r$. If $r>3$, then there is no other elliptic normal
curve containing $\Gamma$.

By contrast, if $r=3$, there are many elliptic normal curves containing
such a $\Gamma$; indeed, the set is parametrized by an open subset of $\P^2$.

\smallskip

\noindent{\sl Proof Sketch.\/} Suppose that $\Gamma\subset E\cap E'$, where
$E$ and $E'$ are elliptic normal curves in $\P^r$, and assume $\Gamma$ is
equivalent to twice the
hyperplane section of $E$.  
It follows that the quadrics vanishing on $E\cup E'$
form a codimension 1 subspace of those vanishing on $E$.  

Suppose $r>3$. The threefold that is the union of 
the secant lines of $E$ has $E$ as its singular locus,
so if $E\neq E'$ there is a secant line to $E$ that is
not secant to $E'$. It follows that we can find distinct
rational normal scrolls $X$ and $X'$, of codimension 2, containing
$E$ and $E'$ respectively.

But the intersection of any two distinct quadrics containing a
codimension 2 scroll is the union of the scroll and a codimension 2
linear subspace; there is no room for another scroll. Therefore $E=
E'$ is the unique elliptic normal curve in $\P^r$ containing
$\Gamma$.\Box

\medskip
We now turn to the classification results. In $\P^1$ the matter
is trivial: every degree 4 scheme is self-associated, and of
course all are arithmetically Gorenstein. The problem is 
already more challenging in $\P^2$ and $\P^3$, and we begin
with some general remarks. Since the 
classification of Gorenstein schemes in these codimensions is
well-known, the difficult point here is to
decide what examples exist that are not arithmetically Gorenstein.

Let $\Gamma\subset\P^r$ be a finite self-associated Gorenstein scheme. 
By \ref{basic-principle} every codegree 2 subscheme of $\Gamma$ spans $\P^r$, 
and in particular $\Gamma$ is nondegenerate. By \ref{vero-gale}
and \ref{Koszu-hom}, $\Gamma$ is Gorenstein unless
 $\Gamma$ is contained in the scheme defined by the 
ideal of minors of a  matrix of the form
$$
\pmatrix{
x_0&\dots&x_t&x_{t+1}&\dots&x_r\cr
0&\dots&0&l_{t+1}&\dots&l_r
},\leqno{(*)}
$$
where $0\le t<r$ and the $l_i$ are linearly independent linear forms.
In particular, $\Gamma_{\red}$ lies in the union of the planes 
$L_1=V(x_0,\dots,x_t)$ and $L_2=V(l_{t+1},\dots,l_r)$. 
If $L_1\cap L_2=\emptyset$, then $\Gamma$ must be 
decomposable, and by \ref{degen-self-assoc}
$\Gamma\cap L_i$ is self-associated in $L_i$ for each $i$.
It seems plausible that something of this sort happens more
generally: 

\problem{} Suppose that $\Gamma$ is self-associated and the ideal of 
$\Gamma$ contains the $2\times 2$ minors of the matrix $(*)$ above. 
Under what circumstances is $\Gamma\cap V(x_0,\dots,x_t)$ self-associated
in its span?

\medskip
From the classification below we see that for $\P^3$, the first
projective space in which a non-trivial example arises, the answer is
``always!'' The following includes a weak result of this type, which
still suffices to eliminate many possibilities:

\lemma{forbidden-cases}
 Suppose that $\Gamma\subset \P^r$ is a self-associated scheme.
\item{a)} If $\Gamma'\subset\Gamma$ has degree $r+d$, then
$\Gamma'$ spans at least a subspace of dimension $d$.
\item{b)} If the homogeneous ideal of $\Gamma$ contains
a product of ideals $(l_1,\dots,l_s)\cdot(m_1,\dots,m_u)$ where 
the $l_i$ are linearly independent linear forms, and 
similarly for the $m_j$, and $s+u>r$, then $2\leq u \leq r-1$
and $2\leq s \leq r-1$.

\proof $a)$ If $\Gamma$ is self-associated
then, as the embedding series is very ample, \ref{bpf} shows
that no subscheme of
$\Gamma$ of degree $2r$ can lie in a hyperplane. Thus no
subscheme of $\Gamma$ of degree $r+d$ can lie in a 
subspace of dimension $d-1$. 

$b)$ Let $\Gamma'=\Gamma\cap V(l_1,\dots,l_s)$. Since the
homogeneous coordinate ring $S_\Gamma$ is Cohen-Macaulay,
any linear form
$l$ vanishing on $\Gamma'$ annihilates the ideal $(m_1,\dots,m_u)$,
so we may harmlessly assume that
the span of $\Gamma'$ is the $(r-s)$-plane
$V(l_1,\dots,l_s)$. The residual
to $\Gamma'$ in $\Gamma$ lies inside $V(m_1,\dots,m_u)$,
so by \ref{basic-principle}, $\Gamma'$ fails by at least
$u$ to impose independent conditions on hyperplanes. If follows
that $\deg(\Gamma')\geq u+(r+1-s)=r+(u-s+1)$. By the result of part $a)$,
we have $r-s \geq u-s+1$, or equivalently 
$u\leq r-1$, one of the desired inequalities. 
By symmetry $s\leq r-1$, and since $s+u>r$, we derive
$2\leq u$ as well.\Box

We can now complete the classification in $\P^2$ and $\P^3$:

\theorem{p2} A finite Gorenstein scheme in $\P^2$ is self-associated
iff it is a complete intersection of 
a conic and cubic.

\proof If $\Gamma$ is arithmetically Gorenstein then since it has codimension 2
it must be a complete intersection. It cannot lie on a line, and it
has degree 6, so it is the complete intersection of a conic and a
cubic. 

If $\Gamma$ is not arithmetically Gorenstein, then the ideal of
$\Gamma$
contains the ideal of the minors of a matrix of the form $(*)$ above. 
In particular it contains $(x_0,\dots,x_t)(l_{t+1},\dots,l_2)$,
so by \ref{forbidden-cases} part $b)$, we get 
$2\leq t+1\leq 1$, a contradiction.
\Box

\proposition{p3} A finite Gorenstein scheme $\Gamma\subset\P^3$
self-associated if and only if either
\item{a)} $\Gamma$ is a complete intersection
of type $(2,2,2)$ (thus the general such $\Gamma$ is a
quadric section of an ``elliptic normal quartic'' curve in $\P^3$), or
\item{b)} $\Gamma$ is cut out by the Pfaffians of a
$5\times 5$ skew symmetric matrix with entries of degrees
$$\pmatrix{
-&-&1&1&1\cr
-&-&1&1&1\cr
1&1&-&2&2\cr
1&1&2&-&2\cr
1&1&2&2&-\cr
},$$
where the dashes denote zero entries, or
\item{c)} There is a smooth quadric $Q$ and a divisor $C$ of
type $(2,0)$ on $Q$ such that $\Gamma$ is a Cartier
divisor on $C$ in the class $4H_C$, where $H_C$ denotes the
hyperplane class (that is,
$\Gamma$ consists of a degree 4 subscheme on each of two disjoint lines,
or a subscheme of degree 8 on a double line meeting the reduced line
in a degree 4 subscheme).

\proof The schemes $\Gamma$ described in $a)$ and $b)$ are 
arithmetically Gorenstein of degree 8, and thus self-associated.
For part $c)$ we may apply \ref{Goppa} to the linearly normal
curve $C$. The restriction map
$\Pic(C) \rTo \Pic(C_{\red})$ is an isomorphism
and writing $H_C$ for the hyperplane class on $C$
it follows that $\Gamma+K_C-H_C=H_C$, so $\Gamma$ is 
indeed self-associated.

For the converse, suppose first that $\Gamma$ is 
arithmetically Gorenstein. By the structure
theorem (Buchsbaum-Eisenbud [1977]) for codimension three arithmetically
Gorenstein schemes, $\Gamma$ has ideal $I_\Gamma$  
generated by the $2n\times 2n$-Pfaffians
of a $(2n+1)\times (2n+1)$ skew symmetric matrix. From the Hilbert function
we know that $I_\Gamma$ contains three quadrics and is moreover 3-regular. 
If $n=1$, the ideal is
generated by these three quadrics, and is thus a complete intersection (case
$a)$). If $n=2$ there must be 2 cubic generators  
in addition to the 3 quadrics, and
the given degree pattern is easy to deduce (case $b)$).  
Finally, if $n>2$, then the
Pfaffians would all
have degree $>2$, which is impossible.

On the other hand, suppose that $\Gamma$
is not arithmetically Gorenstein. By \ref{Koszu-hom}, 
$\Gamma$ lies on the scheme defined by the 
$2\times 2$ minors of a matrix of the form $(*)$. 
By \ref{forbidden-cases} we have $t=1$. 
If $V(l_2,l_3)$ is disjoint from $V(x_0,x_1)$,
then $\Gamma$ lies on the disjoint union of two lines, and is
of course a degree 8 Cartier divisor there. Any two disjoint lines lie
on a smooth quadric, so we are done in this case.

If on the contrary 
$V(l_2,l_3)$ meets or coincides with $V(x_0,x_1)$ then
the matrix $(*)$ can be reduced by
a linear change of variables and columns to the form
$$
\pmatrix{
x_0&x_1&x_{2}&x_3\cr
0&0&x_0&l_3
},
$$
with $l_3$ equal to $x_1, x_2$, or $x_3$.

If $l_3=x_2$ or $l_3=x_3$, then  we see that
$x_0$ corresponds to an element of the socle of the local
ring of $\Gamma$ at the point $V(x_0,x_1,l_3)$, and vanishes
on any component of $\Gamma$ supported away from this point.
Since $\Gamma$ is nondegenerate,
$x_0\neq 0$ in this local ring.
Since $\O_\Gamma$ is Gorenstein, 
$x$ generates the socle of the local ring.
It follows that the line $x=0$ contains a
codegree 1 subscheme of $\Gamma$, contradicting 
part $a)$ of \ref{forbidden-cases}.

Thus we may suppose $l_3=x_1$.
In this case $\Gamma$ lies on a double line on the smooth quadric
$V(x_0x_3-x_1x_2)$, 
and it remains to see that $\Gamma$ is a  Cartier divisor there.
By \ref{forbidden-cases} $a)$ the reduced line can intersect $\Gamma$
in a subscheme of degree at most 4.
Passing to the affine case, we may take a polynomial $f$
in the ideal of $\Gamma$ in the double line which restricts to the reduced
line to define the same scheme of degree 4. Since $f$ is a nonzero-divisor
in the ideal of the double line, it defines a subscheme of degree 8,
and thus $f$ generates the ideal of $\Gamma$ in the double line.
It follows that $\Gamma$ is Cartier, which concludes the proof of the
Proposition. \Box

\remark{free-resolutions} The classification in \ref{p3} is
also the classification by numerical type of the free resolution,
or, as it turns out, by the length of the 2-linear part of the
resolution, the
``resolution Clifford index'' in an obvious sense
(see Eisenbud [1992]). The analogue here of Green's conjecture
might be to show that the
resolution Clifford index is always determined by the
``geometric Clifford index,'' that is, the
types of matrices of the form $(*)$ that arise.
Be this as it may in general, 
the possible free resolutions of $S_\Gamma$ over
$k[x_0,\dots,x_3]$ in cases $a)$, $b)$, and $c)$ respectively,
are
$$
\vbox{\offinterlineskip %\baselineskip=15pt
\halign{\strut\hfil# \ \vrule\quad&# \ &# \ &# \ &# \ &# \ &# \ 
&# \ &# \ &# \ &# \ &# \ &# \ &# \
\cr
degree&\cr
\noalign {\hrule}
0&1 &--&--&--&\cr
1&--&3 &--&--&\cr
2&--&--&3 &--&\cr
3&--&--&--&1 &\cr
\noalign{\bigskip}
%\omit&\multispan{3}{\bf Case a)}\cr
\noalign{\smallskip}
}}
$$
$$
\vbox{\offinterlineskip %\baselineskip=15pt
\halign{\strut\hfil# \ \vrule\quad&# \ &# \ &# \ &# \ &# \ &# \ 
&# \ &# \ &# \ &# \ &# \ &# \ &# \
\cr
degree&\cr
\noalign {\hrule}
0&1 &--&--&--&\cr
1&--&3 &2 &--&\cr
2&--&2 &3 &--&\cr
3&--&--&--&1 &\cr
\noalign{\bigskip}
%\omit&\multispan{3}{\bf Case b)}\cr
\noalign{\smallskip}
}}
$$
$$
\vbox{\offinterlineskip %\baselineskip=15pt
\halign{\strut\hfil# \ \vrule\quad&# \ &# \ &# \ &# \ &# \ &# \ 
&# \ &# \ &# \ &# \ &# \ &# \ &# \
\cr
degree&\cr
\noalign {\hrule}
0&1 &--&--&--&\cr
1&--&4 &4 &1 &\cr
2&--&--&--&--&\cr
3&--&2 &4 &2 &\cr
\noalign{\bigskip}
%\omit&\multispan{3}{\bf Case c)}\cr
\noalign{\smallskip}
}}
$$

\remark{all-are-the-same} Case $b)$ in \ref{p3} above corresponds in fact
to schemes of degree 8 on a (possibly degenerate)
twisted cubic curve.  Indeed, the
determinantal ideal in the first two rows of the $5\times 5$ skew
symmetric matrix must actually have codimension 2 (Proof: Its minors
appear among the 5 minimal generators and are thus linearly
independent.  We may reduce modulo a general linear form and reduce to
a problem in 3 variables.  If the three minors had a common divisor
$x$, then $x$ would be in the socle module of the (reduced) ideal of the
points, which is impossible, as the socle is entirely in degree 3.)
Therefore, the Pfaffians define a scheme of degree 8 on a determinantal curve
of degree 3
in $\P^3$.  Note also, that two general quadrics in the ideal of the
curve define in general an arithmetic genus 1 quartic curve containing
the eight points; but they are not a quadric section of this
quartic.

\smallskip
Here is a geometric description of a special case of case $c)$:

\example{four-doubles}
Suppose $\Gamma\subset\P^3$ is a scheme of degree 8 
consisting of four double points. Suppose further 
that the degree 4 scheme $\Gamma_\red$ is contained in a line
$R$. Then $\Gamma$ is self-associated iff the components
of $\Gamma$ are tangent to four rulings on 
a smooth quadric surface iff the four points of $\P^1$
corresponding to the tangent vectors to $\Gamma$ in the
normal bundle of $R$ have the same cross-ratio as the
corresponding points of $\Gamma_\red$ in $R$.\Box

\smallskip
In $\P^4$ we have a less complete result. The extra hypothesis
of linear general position excludes in particular the non
arithmetically Gorenstein cases such as 
the union of 4 points on a line and 6 points on a conic spanning
a disjoint plane. The result was enunciated by Bath in the reduced
``sufficiently general'' case.

\theorem{bath} Let $\Gamma\subset\P^4$ be a finite, local complete
intersection scheme, which is in linearly general position.
Then $\Gamma$ is self-associated (and in fact arithmetically Gorenstein) 
if and only if either
\item{a)} $\Gamma$ is a quadric section of an elliptic normal quintic curve
(equivalently a hyperplane section of a
non-trigonal canonical curve of genus 6 in $\P^5$), or
\item{b)} $\Gamma$ is a scheme of degree ten on a rational normal
quartic curve.

\proof From the general position hypothesis, 
\ref{lgp-and-omega} and \ref{vero-gale} it follows that
$\Gamma$ is self-associated iff 
$\Gamma$ is arithmetically Gorenstein. In particular, $\Gamma$
fails exactly by one to impose independent condition on
quadrics, and thus $\h^0(\I_\Gamma(2))=6$.

A structure theorem of Kustin-Miller [1985], Herzog-Miller [1986],
and Vasconcelos-Villareal [1986] asserts that a generic
local complete intersection Gorenstein ideal $I$, of grade 4
and deviation 2, is of the type $I=\langle J, f\rangle$,
where $J$ is a (Gorenstein) codimension 3 ideal defined
by the $4\times 4$ Pfaffians of a skew symmetric matrix,
and $f$ is a non-zero divisor modulo $J$.  Thus, in case
the six quadrics in $I_\Gamma$ generate the homogeneous ideal, 
that is $I_\Gamma$ is  Gorenstein  of grade 4 and deviation 2, 
then $\Gamma$ is a quadric section of an
arithmetically Gorenstein scheme $\Lambda\subset\P^4$ defined by
the Pfaffians of a $5\times 5$-skew symmetric matrix
with linear entries, which is case $a)$ in the statement
of the theorem.

Assume now that the quadrics in $\H^0(\I_\Gamma(2))$ do not generate 
the  Gorenstein ideal $I_\Gamma$.
In this case, there are also cubic generators in the ideal, and 
by symmetry their number matches the dimension of 
$\Tor_{3}^S(I_\Gamma,S)_{2}$,
which is thus nonzero. By the ``Strong Castelnuovo Lemma''
of Green [1984], Yanagawa [1994] and Cavaliere-Rossi-Valla [1994];
see also Eisenbud-Popescu [1997], it follows that
$\Gamma$ is divisor of degree 10 on a smooth rational
normal quartic curve, which is case $b)$ in the  statement
of the proposition.\Box

\remark{gen-assoc-p4} Bath [1938] claims that the general
self-associated set in $\P^4$ is a quadric section of a quintic
elliptic normal curve, (case $1)$ in \ref{bath}. (See also Babbage
[1948].) Here is an outline of his argument:

A general self-associated ordered set 
$\Gamma=\{p_1,\dots, p_{10}\}\subset\P^4$
fails by one to impose independent conditions on quadrics, so
$\h^0(\I_\Gamma(2))=6$.
Either $\Gamma$ is contained in a rational normal quartic curve, or
three general quadrics in $\H^0(\I_\Gamma(2))$ meet along a genus 5 canonical
curve $C\subset\P^4$, passing through the set $\Gamma$. In the
latter case, the quadrics in
$\H^0(\I_\Gamma(2))$ cut out a $g^2_6$ residual to $\Gamma$ on the curve $C$.
However, a $g^2_6$ on $C$ is special, and this means that any divisor in
this linear system spans only a $\P^3$. Let  $\Sigma$ be a (general) divisor
in the $g^2_6$, so that $\Sigma$ is reduced, disjoint from $\Gamma$, 
and in linearly
general position in its span (since $C$ is cut out by quadrics, and is not
hyperelliptic). By Castelnuovo's lemma (see for instance \ref{Castel-Gor})
there is a unique twisted cubic curve $D\subset\P^3$ (the linear span of
$\Sigma$) which passes through $\Sigma$.
Now $\Gamma\cup\Sigma$ is the complete intersection of 4 quadrics from
$\H^0(\I_\Gamma(2))$, and since it is only necessary to make a quadric
contain one more point of $D$ for the whole twisted cubic $D$ to
lie on the quadric, it follows that $D$ lies on three independent
quadrics in $\H^0(\I_\Gamma(2))$. They define a complete intersection
curve in $\P^4$, which has as components $D$ and another (arithmetically
Gorenstein) curve $E$ of degree 5, passing through the ten points
$\Gamma=\{p_1,\dots, p_{10}\}$. The curve $E$ is an elliptic quintic curve,
and $\Gamma$ is a  quadric section of it.

In a general self-associated, ordered
set $\Gamma=\{p_1,\dots, p_{10}\}\subset\P^4$, one can always arbitrarily
prescribe the first eight points. As in the proof of \ref{complete},
one sees that among the non-singular
bilinear diagonal forms, for which $\Sigma=\{p_1,\dots, p_{5}\}$ forms
an orthogonal basis, there is a pencil $B_{(s:t)}$ of bilinear forms
for which the points $\Sigma'=\{p_6,p_7,p_8\}$ are also mutually orthogonal.
The conditions that $p_9$ is orthogonal on $\Sigma'$ are expressed by
a system of three bilinear equations in $(s:t)$ and the coordinates of the
ambient
$\P^4$. The system has a solution $B$ iff the $3\times 2$ matrix of the linear
system drops rank iff the point $p_9$ lies on $X\subset\P^4$, the variety
defined by
the maximal minors of the $3\times 2$ matrix (which has linear entries in the
coordinates
of the ambient $\P^4$). For general choices, the bilinear form $B$ is unique
and
nonsingular, and $X\subset\P^4$ is a smooth  rational cubic scroll.
Analogously, the point $p_{10}$ is orthogonal on $\Sigma'$,
with respect to $B$, iff $p_{10}$ lies also on the scroll $X$.
Interpreting $X\subset\P^4$ as the image of $\P^2$ via conics through a point,
one sees readily that there is a pencil of elliptic 
quintic normal curves through
$\{p_1,\dots, p_{8}\}$, which are all bisections for the ruling of the
scroll $X$, and for any given choice of such an elliptic quintic normal
curve $E$ one has to pick $\{p_9, p_{10}\}$ as the intersection points
of $E$ with a ruling of $X$.

\remark{} Mukai [1995] proved that every canonical curve of
genus 7 and Clifford index 3 (i.e., the general canonical curve of genus 7)
is a linear section of the
the spinor variety $S_{10}\subset\P^{15}$ of isotropic
$\P^4$'s in the $8$-dimensional quadric $Q\subset\P^9$.
In the same spirit,
Ranestad and Schreyer [1997] showed that
the ``general empty arithmetically Gorenstein scheme
of degree 12 in $\P^4$'' (i.e., the general graded Artinian Gorenstein
with Hilbert function $(1,5,5,1)$) is always a linear section of
the same spinor variety.  It seems natural to expect a similar result in
our case:

\conjecture{} The general arithmetically Gorenstein, non\-degenerate
zero\--dimensional scheme of degree 12 in $\P^5$ is a linear
section of the spinor variety $S_{10}\subset\P^{15}$.

\remark{} We refer to Babbage [1948] for a description
in the spirit of \ref{gen-assoc-p4}
of the general set of twelve self-associated points in $\P^5$.

\bigskip\bigskip
{
\centerline {\bf References}
\baselineskip=12pt
\parindent=0pt
\frenchspacing
\medskip
%Models:
%\item{} V.~I.~Arnold: A-graded algebras and continued fractions,
%{\sl Communications in Pure and Appl. Math.} {\bf 42} (1989) 993-1000.
%\medskip
%\item{} D.~Cox, J.~Little, D.~O'Shea:
%{\sl Ideals, Varieties and Algorithms},
%Springer, New York, 1992.
%\medskip

%\references

\item{} D.W.~Babbage:
Twelve associated points in $[5]$,
{\sl J. London Math. Soc.} {\bf 23}, (1948), 58--64.
\medskip

\item{} F.~Bath: Ten associated points and quartic and quintic
curves in $[4]$, {\sl J. London Math. Soc.} {\bf 13},
(1938), 198--201.
\medskip

\item{} D.~Bayer, M.~Stillman:
Macaulay: A system for computation in
        algebraic geometry and commutative algebra
Source and object code available for Unix and Macintosh
        computers. Contact the authors, or download from
{\tt ftp:// math.harvard.edu/macaulay} via anonymous ftp.
\medskip

\item{} W.~Bruns, J.~Herzog: Cohen-Macaulay rings. 
Cambridge Studies in Advanced Mathematics, {\bf 39}. 
Cambridge University Press, Cambridge, 1993. 
\medskip

%\item{} D.~Buchsbaum, D.~Eisenbud:
%Some structure theorems for finite free resolutions, {\sl
%Advances in Math.} {\bf 12}, (1974), 84--139.
%\medskip

\item{} D.~Buchsbaum, D.~Eisenbud:
Algebra structures for finite free resolutions, and some structure
theorems for ideals of codimension $3$.
{\sl Amer. J. Math.} {\bf 99}, (1977), no. 3, 447--485.
\medskip

\item{} M.P.~Cavaliere, M.E.~Rossi, G.~Valla:
The Strong Castelnuovo Lemma for zerodimensional schemes,
in ``Zero-dimensional schemes (Ravello, 1992)'', 53--63,
de Gruyter, Berlin, 1994.
\medskip

\item{} M.P.~Cavaliere, M.E.~Rossi, G.~Valla:
Quadrics through a set of points and their syzygies,
{\sl Math. Z.}, {\bf 218}, (1995), 25--42.
\medskip

\item{} G.~Castelnuovo:
Su certi gruppi associati di punti,
{\sl Ren. di Circ. Matem. Palermo} {\bf 3},
(1889), 179--192.
\medskip

\item{} A.B.~Coble: Point sets and allied Cremona groups I,
{\sl Trans. Amer. Math. Soc.} {\bf 16}, (1915) 155--198.
\medskip

\item{} A.B.~Coble: Point sets and allied Cremona groups II,
{\sl Trans. Amer. Math. Soc.} {\bf 17}, (1916) 345--385.
\medskip

\item{} A.B.~Coble: Point sets and allied Cremona groups III,
{\sl Trans. Amer. Math. Soc.} {\bf 18}, (1917) 331-372.
\medskip

\item{} A.B.~Coble: Associated sets of points,
{\sl Trans. Amer. Math. Soc.} {\bf 24}, (1922) 1--20.
\medskip

\item{} A.B.~Coble:
{\sl Algebraic Geometry and Theta Functions},
American Mathematical Society, New York, 1929.
\medskip

\item{} J.R.~Conner: Basic systems of rational norm-curves,
{\sl Amer J. Math.} {\bf 32}, (1911) 115--176.
\medskip

\item{} J.R.~Conner: The Rational Sextic Curve, and the
Caylay Symmetroid, {\sl Amer J. Math.} {\bf 37}, (1915) 28--42.
\medskip

\item{} F.~Davide: The Gale transform for sets of points on the
Veronese surface, Preprint 1997
\medskip

\item{} E.~D.~Davis, A.~V.~Geramita, F.~Orecchia: Gorenstein
algebras and the Cayley-Bacharach theorem,
{\sl Proc. Amer. Math. Soc.}, {\bf 93}, (1985), 593--597.
\medskip

\item{} H.~Dobriner: \"Uber das Ra\"umliche Achteck welches
die Schnittpunkte dreier Oberfl\"achen zweiter Ordnung bilden,
{\sl Acta Math.} {\bf 12}, (1889), 339--361.
\medskip

\item{} I.~Dolgachev, M.~Kapranov:
Arrangements of hyperplanes and vector bundles on 
$\P^n$, {\sl Duke Math. J.}, {\bf  71},
(1993), no. 3, 633--664.
\medskip

\item{} I.~Dolgachev, D.~Ortland: Points sets in
projective spaces and theta functions, {\sl Ast\'erisque}
{\bf 165}, (1988).
\medskip

\item{} J.~Edmonds: Minimum partition of a matroid into
independent subsets, {\sl J. of Research of the
national Bureau of Standards - B. Mathematics and
Mathematical Physics}, {\bf 69} B, (1965), 67--72.
\medskip

\item{} S.~Ehbauer:
Syzygies of points in projective space and applications,
in ``Zero-dimensional schemes (Ravello, 1992)'', 145--170,
de Gruyter, Berlin, 1994.
\medskip

\item{} D.~Eisenbud:
Linear sections of determinantal varieties,
{\sl Amer. J. Math.} {\bf 110} (1988), no. 3, 541--575.
\medskip

\item{} D.~Eisenbud: Green's conjecture: an orientation for
algebraists, in ``Free Resolutions in Commutative Algebra and
Algebraic Geometry'', {\sl Research Notes in Mathematics}, ({\bf 2}),
51--78, Jones and Bartlett Publishers, Boston, 1992.
\medskip

\item{} D.~Eisenbud:
{\sl Commutative Algebra with a View Toward Algebraic Geometry},
Springer, New York, 1995.
\medskip

\item{} D.~Eisenbud, M.~Green, J.~Harris:
Cayley-Bacharach theorems and conjectures.
{\sl Bull. Am. Math. Soc}, {\bf 33}, (1996), no. 3, 295--324.
\medskip

\item{} D.~Eisenbud, M.~Green, J.~Harris:
Higher Castelnuovo theory, in ``Journ\'es
de G\'eom\'etrie Alg\'ebrique d'Orsay (Orsay, 1992)''.
{\sl Ast\'erisque} {\bf 218} (1993), 187--202.
\medskip

\item{} D.~Eisenbud, J.~Harris:
Finite projective schemes in linearly general position.
{\sl J. Algebraic Geom.}, {\bf 1} (1992), no. 1, 15--30.
\medskip

\item{} D.~Eisenbud, J.~Harris:
 An intersection bound for rank $1$ loci, with
applications to Castelnuovo and Clifford theory.
{\sl J. Algebraic Geom.} {\bf 1} (1992), no. 1, 31--59.
\medskip

\item{} D.~Eisenbud, J-H.~Koh: Remarks on points in
a projective space, in
{\sl Commutative Algebra (Berkeley, CA, 1987)},
{\sl Math. Sci. Res. Inst. Publ.}, {\bf 15}, 157--172,
Springer 1989
\medskip

\item{} D.~Eisenbud, J-H.~Koh, M.~Stillman: Determinantal
equations for curves of high degree, {\sl Amer. J. of Math.}, 
{\bf 110}, (1988), 513--539.
\medskip

\item{} D.~Eisenbud, S.~Popescu:
Gale Duality and Free Resolutions of Ideals of Points,
Preprint  1996, (alg-geom/9606018).
\medskip

\item{} D.~Eisenbud, S.~Popescu:  Syzygy Ideals for Determinantal Ideals
and the Syzygetic Castelnuovo Lemma.  Preprint, 1997.
\medskip

\item{} F.~Flamini: Inductive construction of self-associated sets of
points, Preprint 1997.
\medskip

\item{} D.~Gale: Neighboring vertices on a convex polyhedron,
in ``Linear Inequalities and Related Systems'' (H.W.~Kuhn
and A.W.~Tucker, eds.), {\sl Annals of Math. Studies} {\bf 38},
255--263, Princeton Univ. Press, 1956.
\medskip

\item{} A.~Gimigliano: On Veronesean surfaces, {\sl Indag. Math.},
A {\bf 92}, (1), (1989), 71--85.
\medskip

\item{} V.D.~Goppa: A new class of linear error-correcting codes,
{\sl Problems of Information Transmission} {\bf 6}, (1970) 207--212.
\medskip

\item{} V.D.~Goppa: Codes and Information,
{\sl Russian Math. Surveys} {\bf 39}, (1984) 87--141.
\medskip

\item{} D.~Grayson, M.~Stillman:
Macaulay 2: a software system devoted to supporting research
in algebraic geometry and commutative algebra.
Contact the authors, or download source and binaries from
{\tt ftp://ftp.math.uiuc.edu/Macaulay2 } via anonymous ftp.
\medskip

\item{} M.~Green: Koszul cohomology and the geometry of projective varieties I,
{\sl J. Differential Geom.} {\bf 19} (1984), no. 1, 125--171
\medskip

\item{} M.~Green, R.~Lazarsfeld:
Some results on the syzygies of finite sets and algebraic curves,
{\sl Compositio Math.} {\bf 67}, (1988), no. 3, 301--314.
\medskip

\item{} L.~Gruson, Ch.~Peskine: 
Courbes de l'espace projectif: vari\'et\'es de s\'ecantes,
in {\sl Enumerative geometry and classical algebraic geometry (Nice, 1981)},
Progr. Math., {\bf 24}, Birkh\"auser (1982), 1--31.
\medskip

\item{} J.~Herzog, M.~Miller: Gorenstein ideals of deviation two, {\sl Comm.
Algebra}, {\bf 13}, (1985), no. 9, 1977--1990.
\medskip

\item{} L.O.~Hesse: De octo punctis intersectionis trium
superficierum secundi ordinis, Dissertatio, (1840), Regiomonti.
\medskip

\item{} L.O.~Hesse: De curvis et superficiebus secundi ordinis,
{\sl J. reine und angew. Math.} {\bf 20}, (1840), 285--308.
\medskip

\item{} L.O.~Hesse: Ueber die lineare Construction des achten Schnittpunktes
dreier Oberfl\"achen zweiter Ordnung, wenn sieben Schnittpunkte
derselben gegeben sind,
{\sl J. reine und angew. Math.} {\bf 26}, (1840), 147--154.
\medskip

\item{} K.~Hulek, Ch.~Okonek, A.~van~de~Ven: Multiplicity-2 
structures on Castelnuovo surfaces, {\sl Ann. Scuola Norm. Pisa} (4),
{\bf 13}, (1986), no. 3, 427--448.
\medskip

\item{} Igusa: On Siegel modular forms of genus two, I and II {\sl Amer. J.
Math.}, {\bf 86}, (1964), 219--246, and  {\bf 88}, (1964), 392--412.
\medskip

\item{} M.~Kapranov:
Chow quotients of Grassmannians. I,  in I. M. Gel'fand Seminar,
{\sl Adv. Soviet Math.}, {\bf 16}, Part 2,
Amer. Math. Soc., Providence, RI, (1993), 29--110.
\medskip

\item{} M.~Kreuzer: Some applications of the canonical module of
a $0$-dimensional scheme, in 
{\sl ``Zero-dimensional schemes'' (Ravello, 1992)},
243--252, de Gruyter, Berlin, 1994.
\medskip

\item{} M.~Kreuzer: On $0$-dimensional complete intersections,
{\sl Math. Ann.}, {\bf 292}, (1992), no. 1, 43--58.
\medskip

\item{} A.R.~Kustin, M.~Miller: Classification of the {\rm Tor}-algebras of
codimension four Gorenstein local rings, {\sl Math. Z.}, {\bf 190}, (1985),
341--355.
\medskip

\item{} J.M.~Landsberg: On an unusual conjecture of Kontsevich
and variants of Castelnuovo's lemma, Preprint 1996, (alg-geom/9604023).
\medskip

\item{} J.H.~van Lint, G.~van der Geer:
{\sl Introduction to Coding Theory and Algebraic Geometry},
DMV Seminar, Band {\bf 12}, Birkh\"auser, Basel, 1988.
\medskip

\item{} F.~Meyer (edited by):
Encyclop\'edie des Sciences Math\'ematiques Pures et Appliqu\'ees,
publi\'ee sous les auspices des acad\'emies des sciences de
G\"ottingue, de Leipzig, de Munich et de Vienne avec la
collaboration de nombreux savants. (\'edition fran\c aise,
r\'edig\'ee et publi\'ee d'apr\'es l'\'edition allemande
sous la direction de J. Molk), Gauthier-Villars, Paris,
\& Teubner, Leipzig,
reprinted by Jacques Gabais, Paris (1992))
\medskip

\item{} S.~Mukai: Curves and Symmetric Spaces, I, {\sl Amer. J. Math.},
{\bf 117}, (1995), 1627--1644.
\medskip

\item{} B.~Ramamurti: On ten associated points in $[4]$,
{\sl Proc Indian Acad. Sci.}, Sect. A, {\bf 16}, (1942), 191-192.
\medskip

\item{} K.~Ranestad, F.-O.~Schreyer: Varieties of sums of powers,
preprint 1997, (alg-geom/9801110).
\medskip

\item{} J.~Rathmann: Double structures on Bordiga surfaces, 
{\sl Comm. Algebra}, {\bf 17}, (1989), no. 10, 2363--2391. 
\medskip

\item{} T.G.~Room: {\sl The Geometry of Determinantal Loci},
Cambridge Univ. Press, London, 1938.
\medskip

\item{} J.~Rosanes: \"Uber linear abh\"angige Punktsysteme,
{\sl J. reine und angew. Math.} {\bf 88}, (1880), 241--272.
\medskip

\item{} J.~Rosanes: Zur Theorie der reciproken Verwandschaften,
{\sl J. f\"ur die Reine und Angew. Math.} {\bf 90}, (1881), 303--321.
\medskip

\item{} V.~Shokurov: The Noether-Enriques theorem on canonical
curves, {\sl Mat. Sbornik} {\bf 15}, (1971), 361--403 (engl. transl).
\medskip

\item{} K.G.~Chr.~von~Staudt:
{\sl Beitr\"age zur geometrie der Lage}, fasc. {\bf 3}, (1860), p.373.
\medskip

\item{} J.~Stevens: On the number of points determining a canonical curve,
{\sl Nederl. Akad. Wetensch. Indag. Math.}, {\bf 51}, (1989), no. 4,
485--494.
\medskip

\item{} D.~Struik,: {\sl A source book in mathematics,  1200--1800},
Harvard Univ. Press, Cambridge, MA, 1969.
\medskip

\item{} R.~Sturm: Das Problem der Projektivit\"at und seine Andwendungen
auf die Fl\"achen zweiten Grades, {\sl Math. Ann} {\bf 1}, (1869),
533--574.
\medskip

\item{} R.~Sturm: \"Uber Collineation und Correlation,
{\sl Math. Ann} {\bf 12}, (1877), 254--368.
\medskip

\item{} V.W.~Vasconcelos, R.~Villarreal: On Gorenstein ideals
of codimension four, {\sl Proc. Amer. Math. Soc.}, {\bf 98}, (1986),
no. 2, 205--210.
\medskip

\item{} T.~Weddle: {\sl Cambr. Dublin Math. J.} {\bf 5}, (1850), 58--234.
\medskip

\item{} H.~Whitney: Some combinatorial properties of complexes,
{\sl Proc. Nat. Acad. Sci. U.S.A.} {\bf 26}, (1940), 143--148.
\medskip

\item{} K.~Yanagawa: Some generalizations of Castelnuovo's lemma
on zero-dimensional schemes. {\sl J. Alg.} {\bf 170} (1994)
429--431.
\medskip

\item{} H.G.~Zeuthen: Note sur les huit points d'intersection
de trois surfaces du second ordre, {\sl Acta Math.} {\bf 12},
(1889), 362--366.
\medskip

}
\bigskip
\vbox{\noindent Author Addresses:
\smallskip
\noindent{David Eisenbud}\par
\noindent{Department of Mathematics, University of California, Berkeley,
Berkeley CA 94720}\par
\noindent{eisenbud@math.berkeley.edu}
\smallskip
\noindent{Sorin Popescu}\par
\noindent{Department of Mathematics, Columbia  University,
New York, NY 10027}\par
\noindent{psorin@math.columbia.edu}\par
}

\bye